\documentclass[12pt]{amsart}

\usepackage{environ, xcolor}
\NewEnviron{commentA}{}

\newcommand\Aoff{\RenewEnviron{commentA}{}}
\Aoff{} %Change this to Aoff for professional version, Aon to see extra comments

\usepackage{amsmath, amssymb}
\usepackage{array}
\usepackage[frame,cmtip,arrow,matrix,line,graph,curve]{xy}
\usepackage{graphpap, color, paralist, pstricks}
\usepackage[mathscr]{eucal}
\usepackage[pdftex]{graphicx}
\usepackage[pdftex,colorlinks,backref=page,citecolor=blue]{hyperref}
\usepackage{cleveref}
\usepackage{tikz, verbatim}

\newtheorem{theorem}{Theorem}[section]
\newtheorem{proposition}[theorem]{Proposition}
\newtheorem{corollary}[theorem]{Corollary}
\newtheorem{lemma}[theorem]{Lemma}

\theoremstyle{definition}
\newtheorem{definition}[theorem]{Definition}
\newtheorem{example}[theorem]{Example}

\newtheorem{question}[theorem]{Question}
\newtheorem{remark}[theorem]{Remark}

\newtheorem{setup}[theorem]{Setup}

\usepackage[margin=1in]{geometry} 
\usepackage{amsmath,amssymb,mathtools, mathrsfs,esint,tikz-cd,verbatim}
\newcommand{\R}{{\mathbb R}}

\newcommand{\Z}{{\mathbb Z}}

\newcommand{\F}{{\mathbb F}}
\newcommand{\C}{{\mathbb C}}

\newcommand{\on}[1]{\operatorname{#1}}

\newcommand{\Pic}{{\on{Pic}}}

\DeclareFontFamily{U}{mathb}{\hyphenchar\font45}
\DeclareFontShape{U}{mathb}{m}{n}{
      <5> <6> <7> <8> <9> <10> gen * mathb
      <10.95> mathb10 <12> <14.4> <17.28> <20.74> <24.88> mathb12
      }{}
\DeclareSymbolFont{mathb}{U}{mathb}{m}{n}
\DeclareMathSymbol{\righttoleftarrow}{3}{mathb}{"FD}


\title{Hypersurfaces with Large Automorphism Groups}
\author{Louis Esser}
\author{Jennifer Li}

\address{Department of Mathematics, Princeton University, Fine Hall, Washington Road, Princeton, NJ 08544-1000, USA}
\email{esserl@math.princeton.edu}

\email{jenniferli@princeton.edu}

\begin{document}

\begin{abstract}
We find sharp upper bounds on the order of the automorphism group of a 
hypersurface in complex projective space in every dimension and degree.
In each case, we prove that the hypersurface realizing the upper bound
is unique up to isomorphism and provide explicit generators for the
automorphism group.
\end{abstract}

\maketitle

In this paper, we consider the following question for any positive 
integers $n,d$:

\begin{question}
    \label{mainq}
What is the largest possible automorphism group of a smooth
hypersurface $X$ of degree $d$ in complex projective space $\mathbb{P}_{\C}^{n+1}$?
\end{question}

When $d$ is $1$ or $2$, all the smooth hypersurfaces for an arbitary 
fixed $n$ are isomorphic to one another 
and furthermore have infinite automorphism group.  \Cref{mainq}
is therefore meaningful only when $d \geq 3$.
When $(n,d) = (1,3)$, $X$ is a smooth genus $1$ curve and has an 
infinite automorphism group.  When $(n,d) = (2,4)$, $X$ is a 
quartic K3 surface, and some of these also have infinite
automorphism group (see, e.g. \cite[Proof of Theorem 4]{MM}).

In all the remaining cases, $\mathrm{Aut}(X)$ is finite
by results of Matsumura-Monsky
\cite{MM}.

For almost all values of $n,d$, the precise upper bound for
$|\mathrm{Aut}(X)|$ was not known.  Our main result
is to calculate the upper bound exactly for every $n$ and $d$, 
completely answering \Cref{mainq}.

\begin{theorem}
\label{mainthm}
Let $X_d \subset \mathbb{P}^{n+1}_{\C}$ be a smooth hypersurface
of degree $d$, where $d \geq 3$ and $(n,d) \neq (1,3), (2,4)$.
Then $|\mathrm{Aut}(X)| \leq (n+2)!d^{n+1}$ except for the following
six cases, where the indicated bound holds:
\begin{enumerate}
    \item $n = 1, d = 4$: $|\mathrm{Aut}(X)| \leq 168$,
    \item $n = 1, d = 6$: $|\mathrm{Aut}(X)| \leq 360$,
    \item $n = 2, d = 6$: $|\mathrm{Aut}(X)| \leq 6912$,
    \item $n = 2, d = 12$: $|\mathrm{Aut(X)}| \leq 86400$,
    \item $n = 4, d = 6$: $|\mathrm{Aut}(X)| \leq 6531840$,
    \item $n = 4, d = 12$: $|\mathrm{Aut}(X)| \leq 186624000$.
\end{enumerate}
Each of these bounds is sharp, and the hypersurface $X$ achieving
the bound is unique up to isomorphism for every pair $(n,d)$.
\end{theorem}

The generic bound is achieved by the Fermat hypersurface.  A stronger
and more precise form of \Cref{mainthm} is given in \Cref{mainthmfull},
where we explicitly describe the hypersurfaces and groups realizing
each upper bound.

Right before this paper appeared, another work was announced by 
S.~Yang, X.~Yu, and Z.~Zhu \cite{YangYuZhu}
with a different proof of \Cref{mainthm}.
Both results use Collins' works on the classification of subgroups
of $\mathrm{GL}_N(\C)$ \cite{Collinsprim,Collins}, but the 
methods otherwise differ.

\textit{Outline.} The outline of this paper is as follows.  
In \Cref{sect:preliminaries},
we introduce the necessary preliminaries in 
the representation theory of finite groups.  
\Cref{sect:mainresults} summarizes
what was previously known and gives
a more precise statement of \Cref{mainthm}, including
equations and generators for the hypersurfaces
with largest automorphism group in the exceptional cases.
\Cref{sect:bounds} introduces the machinery
to prove the upper bounds in \Cref{mainthm} and completes the
proof in most dimensions.  \Cref{sect:proof} finishes the proof 
in the remaining cases.

\textit{Acknowledgements.} The authors would like to thank J\'{a}nos
Koll\'{a}r for helpful discussions and comments.

\section{Preliminaries}
\label{sect:preliminaries}

This section summarizes the material from group and representation
theory that will be required for the proof of 
\Cref{mainthm}.  Throughout, we will work over the complex
numbers $\C$, so all representations of finite groups will
be complex.  All representations in the paper will also be 
finite-dimensional.

\subsection{Group Theory and Notations}

We'll use the following standard notations for various groups:
\begin{itemize}
    \item $\mu_r$, the cyclic group of order $r$,
    \item $r \coloneqq \mu_r$, where no confusion arises,
    \item $S_r$, the symmetric group on $r$ elements,
    \item $A_r$, the alternating group on $r$ elements,
    \item $H_{p^m}$, the Heisenberg group of unipotent 
    $3 \times 3$ matrices over $\F_p$ when $m = 1$, or
    more generally an extraspecial group of order $p^{1+2m}$,
    \item $\mathrm{GL}_n(\F)$, the general linear group 
    of invertible $n \times n$ matrices over the field $\F$,
    \item $\mathrm{SL}_n(\F)$, the special linear group
    of $n \times n$ matrices over $\F$ of unit determinant,
    \item $\mathrm{PGL}_n(\F)$, the general projective
    linear group $\mathrm{GL}_n(\F)/\F^{\times}$,
    \item $\mathrm{Sp}_{2n}(\F)$, the symplectic group of 
    $2n \times 2n$ matrices over $\F$,
    preserving a non-degenerate skew-symmetric bilinear form,
    \item $\mathrm{PSp}_{2n}(\F)$, the projective symplectic
    group of $2n \times 2n$ matrices over $\F$.
\end{itemize}

The groups $\mathrm{SL}_n(\F)$, etc., are written
$\mathrm{SL}_n(q)$, etc., when $\F = \F_q$ is the finite field with 
$q$ elements.

There are a few other groups that appear in the paper. 
In particular, the simple groups $\mathrm{PSU}_4(3)$ and $O_8^+(2)$ 
of Lie type appear as Jordan-H\"{o}lder
constituents of important subgroups of $\mathrm{GL}_6(\C)$
and $\mathrm{GL}_8(\C)$, respectively.  Since there are simple
descriptions for the matrix generators for these subgroups,
we don't provide additional definitions related to finite groups
of Lie type, and instead refer to the references, e.g.
\cite{Collinsprim}.

\begin{definition} \label{centralExtension_def}
We say that $G$ is an \textit{extension} of $\bar{G}$ by $H$
if there exists an exact sequence
\begin{center}
$1 \rightarrow H \rightarrow G \rightarrow \bar{G} \rightarrow 1$.
\end{center}
If $H$ is contained in $Z(G)$, then the extension is \textit{central}.
\end{definition}

We often use the notation from the group ATLAS to write down finite groups;
for instance, the notation $3.A_6$ refers to the perfect extension of 
$A_6$ by the cyclic group $\mu_3$, also known as the Valentiner group.

Suppose that $\rho: G \rightarrow \mathrm{GL}(V)$ is a 
representation of the finite group $G$. We say that the 
representation $\rho$ is \textit{irreducible} if it does not
decompose nontrivially as a direct sum of representations. 
We frequently make use of the following consequence of Schur's lemma:

\begin{lemma}
\label{irrep_center}
Suppose $g \in Z(G)$ is an element of the center of $G$
and $\rho: G \rightarrow \mathrm{GL}(V)$ is an irreducible
representation.  Then $\rho(g)$ is a scalar matrix.
\end{lemma}

In particular, the center of an irreducible $G \subset \mathrm{GL}(V)$
coincides with its intersection with the scalar matrices.

We say that 
the representation $\rho$ as above is \textit{imprimitive} if there
is a nontrivial direct sum decomposition $V = V_1 \oplus \cdots \oplus V_r$
such that the action of $G$ permutes $V_1,\ldots,V_r$.
That is, for each $g \in G$ and each $i = 1,\ldots,r$,
$g(V_i) = V_j$ for some $j$.  (Note that for our purposes, reducible
representations are imprimitive.)

A group representation $\rho$ is called \textit{primitive}
if it is not imprimitive.  We call a subgroup $G$ of $\mathrm{GL}(V)$
primitive when the corresponding representation is so.

Given a representation $\rho: G \rightarrow \mathrm{GL}(V)$, there is an 
induced action of $G$ on the polynomial algebra over $V$ by
$(g \cdot f)(\textbf{x}) \coloneqq f(\rho(g)^{-1} \textbf{x})$.
A polynomial $f$ is an 
\textit{invariant} of the representation 
$\rho$ if $g \cdot f = f$ for each $g \in G$.
We say $f$ is a \textit{semi-invariant} of $\rho$ 
if $g \cdot f = \chi(g)f$ for 
some linear character $\chi$ of $G$.

\begin{lemma}
\label{semi-invariant}
Let $\rho: G \rightarrow \mathrm{GL}(V)$ be a representation. Any semi-invariant polynomial of $\rho$
is an invariant polynomial of $\rho|_N$, where $N \subset G$
is a normal subgroup with $G/N$ abelian. In particular, 
if the group $G$ has no nontrivial abelian quotients, i.e., is 
perfect, then 
every semi-invariant polynomial of a representation of $G$ is invariant.

Conversely, if $N \subset G$ is a normal subgroup with $G/N$ finite abelian
and there is an invariant polynomial of degree $d$ for $\rho|_N$, there is a
semi-invariant polynomial of degree $d$ for $\rho$.
\end{lemma}

\begin{proof}
Any semi-invariant polynomial $f$ satisfies $g \cdot f = \chi(g)f$, where
$\chi$ is some one-dimensional representation of $G$.  The polynomial
$f$ is therefore invariant for the restriction of $\rho$ to the kernel $N$
of $\chi$.  Since $G/N \subset \C^*$, it is abelian.

For the last statement, suppose $f$ of degree $d$ is invariant under $\rho|_N$.  Then the span of the orbit of $f$ under the $G$-action
is a finite-dimensional vector space of polynomials of degree $d$
with induced action by $G/N$. Since $G/N$ is abelian, there is an eigenvector
for this action, which is a semi-invariant of degree $d$ for $\rho$.
\end{proof}

\subsection{Primitive Subgroups of the General Linear Group}

A result of Collins bounds the size
of primitive subgroups of $\mathrm{GL}_N(\C)$ for every $N$.
The proof of this result utilizes, at least in a small way,
the classification of finite simple groups.

\begin{theorem}{{{\cite[Theorem A]{Collinsprim}}}}
\label{primgroups}
Let $G \subset \mathrm{GL}_N(\C)$ be a primitive subgroup
for $N \geq 2$.
Then the index $[G : Z(G)]$ of the center of $G$ is 
bounded for fixed $N$.  Moreover, if $[G : Z(G)]$ achieves
the upper bound, $G/Z(G) \cong S_{N+1}$, unless $N = 2,3,4,5,6,7,8,9,$
or $12$.  In these exceptional cases, the bounds on $[G : Z(G)]$ 
are listed in \Cref{primtable}.  When
the bound is achieved, $G/Z(G) \cong H/Z(H)$, where $H$
is the group shown in the third column.
\end{theorem}

\begin{table}[]
    \centering
    \begin{tabular}{c|c|c|c}
       $N$ & Largest $[G : Z(G)]$ & $H$   \\
       \hline
       $2$ & $60$ & $2.A_5$ \\ 
       $3$ & $360$ & $3.A_6$ \\
       $4$ & $25920$ & $\mathrm{Sp}_4(3)$ \\
       $5$ & $25920$ & $\mathrm{PSp}_4(3)$ \\
       $6$ & $6531840$ & $6_1.\mathrm{PSU}_4(3).2_2$ \\ 
       $7$ & $1451520$ & $\mathrm{Sp}_6(2)$ \\
       $8$ & $348364800$ & $2.O^+_8(2).2$ \\
       $9$ & $4199040$ & $3^{1+4}.\mathrm{Sp}_4(3)$ \\
       $12$ & $448345497600$ & $6.\mathrm{Suz}$
    \end{tabular}
    \caption{The exceptional primitive subgroups of 
    $\mathrm{GL}_N(C)$.  The second column gives the upper bound
    on $[G: Z(G)]$, while the third column gives the normal structure
    of a group $H$ with minimal $|Z(H)|$.}
    \label{primtable}
\end{table}

\begin{remark}
The corresponding table in \cite{Collinsprim} has a typo in the row for
$N = 8$, which is corrected in \Cref{primtable}: the correct value for
$[H:Z(H)]$ for $H = 2.O_8^+(2).2$ is $348364800$, not $348368800$.
\end{remark}

\begin{definition}
\label{Xifunction}
Let $\Xi(N)$ denote the function whose output
is the upper bound on $[G:Z(G)]$ over all primitive subgroups $G$ 
of $\mathrm{GL}_N(\C)$.  Then $\Xi(1) = 1$ and $\Xi(N) = (N+1)!$
for each $N \geq 2$ apart from the exceptional values 
of $N$ from \Cref{primtable}, at which $\Xi$ takes the corresponding
value from the second column.
\end{definition}

We also need some more specific information about the classification
of primitive subgroups of $\mathrm{GL}_N(\C)$ for certain small
values of $N$.  Such groups have been completely classified for $N \leq 11$
by the work of many authors (see \cite[Section 1]{Collinsprim} for some 
history).  These results are
frequently phrased as classifications of \textit{unimodular}
subgroups (i.e., subgroups of $\mathrm{SL}_N(\C)$).  This distinction
is relatively unimportant for us, by the following lemma:

\begin{lemma}
\label{GLSL}
Let $G \subset \mathrm{GL}(V)$ be a finite primitive subgroup.
There is a finite primitive group $H \subset \mathrm{SL}(V)$ 
such that $G/Z(G) \cong H/Z(H)$ (so in particular, $G$ and $H$ have
the same image in $\mathrm{PGL}(V)$) and any polynomial
semi-invariant for $G$ is a polynomial semi-invariant for $H$.
\end{lemma}

\begin{proof}
We construct $H$ as follows.  Let $g_1,\ldots,g_s$ be a minimal generating
set for $G$.  For each $i = 1,\ldots,s$, choose a scalar $a_i \in \C^*$
such that $\det(a_ig_i) = 1$.  Then define
$H \coloneqq \langle a_1 g_1, \ldots, a_s g_s \rangle \subset 
\mathrm{SL}(V)$. The group $H$ is primitive (and in particular irreducible),
because if there were a decomposition $V = V_1 \oplus \cdots \oplus V_r$
with $V_1,\ldots,V_r \subset V$ permuted by $H$, then they would be permuted
by each $a_i g_i$ and hence every $g_i$, contradicting the assumption that $G$ is primitive.

The centers $Z(G)$ and $Z(H)$ coincide with the scalar subgroups of $G$
and $H$ by \Cref{irrep_center}.  Since the images of $g_i$ and $a_i g_i$
are the same in $\mathrm{PGL}(V)$, the images $G/Z(G)$ and $H/Z(H)$
of $G$ and $H$, respectively, in this group are isomorphic.  Finally,
if $f$ is a polynomial semi-invariant of $G$, it follows from the definition
of $H$ that $h \cdot f$ is a scalar multiple of $f$ for any $h \in H$.
\end{proof}

In light of \Cref{GLSL}, we can replace any primitive subgroup of 
$\mathrm{GL}_N(\C)$ with a unimodular subgroup having very similar 
properties.  Since this in general replaces invariants with
semi-invariants, we phrase all the invariant theory results below in terms of 
semi-invariant polynomials.

Also, we note that some of the references (e.g., \cite{Feitdim8})
actually classify \textit{quasiprimitive} 
subgroups of $\mathrm{GL}_N(\C)$, namely those for which the restriction
of the representation to any normal subgroup is a direct sum of copies
of exactly one irreducible representation of the subgroup.  This is a 
weaker property than primitive (see, e.g., \cite{Wales} for more),
so the primitive groups are certainly included in the list. From now on,
we won't mention quasiprimitivity further; at worst, this could result
in some unrealized cases in the lists.

The proof of \Cref{mainthm} requires additional information about primitive
subgroups of $\mathrm{GL}_N(\C)$ for $N = 2,3,4,5,6,$ and $8$.
All these cases besides $N = 8$ are summarized in \cite[Section 8.5]{Feit}.
For $N = 2, 3, 4$, we also refer to the lists of primitive subgroups
with explicit matrix generators in \cite{Blichfeldt}.
The degrees of (semi-)invariant polynomials may be computed in \texttt{MAGMA} or
\texttt{GAP}.\footnote{\texttt{MAGMA} can be used to compute (semi-)invariant polynomials
given a set of matrix generators for a finite group. \texttt{GAP} has
useful databases of finite groups and their representations, which can
be used to compute Molien series, and hence degrees of 
(semi-)invariant polynomials.  Examples of the code used to perform these computations
can be found on the first author's website.}
(Semi-)invariant polynomials exist for any finite group representation;
for instance, taking the product of a $G$-orbit of linear polynomials
shows that there exists such an invariant of degree $|G|$.
The next several propositions list the properties we need about
primitive representations and their invariant theory.

\begin{proposition}
\label{primitive2}
Suppose $G$ is a primitive subgroup of $\mathrm{GL}_2(\C)$.
Then $G/Z(G) \cong H/Z(H)$, where $H$ is one of the following:
\begin{enumerate}
    \item \label{prim2:A5}
    $H = 2.A_5$, the binary icosahedral group, with $[H:Z(H)] = 60$;
    \item \label{prim2:S4}
    $H = 2.S_4$, the binary octahedral group, with $[H:Z(H)] = 24$;
    \item \label{prim2:A4}
    $H = 2.A_4$, the binary tetrahedral group, with $[H:Z(H)] = 12$.
\end{enumerate}
The smallest possible degree of a semi-invariant polynomial is $12$ in case
\eqref{prim2:A5}, $6$ in case \eqref{prim2:S4}, and $4$ in case 
\eqref{prim2:A4}.
In appropriate coordinates,
the unique degree $12$ semi-invariant polynomial of $H = 2.A_5$
may be written
$$x y (x^{10} + 11x^5 y^5 - y^{10}),$$
and the unique degree $6$ semi-invariant polynomial of $H = 2.S_4$ may be
written
$$x^5 y - x y^5.$$
\end{proposition}

\begin{proof}
The classification of subgroups of $\mathrm{GL}_2(\C)$ is well-known
and dates back to work of Klein in the 19th century \cite{Klein2}.
The semi-invariant polynomials for the various groups are discussed in
\cite[Section 0.13]{PV} and \cite[Section 1.2]{Mukai}.  Their degrees
correspond to geometric features of the corresponding polyhedron.  In
particular, the number of vertices of the icosahedron, octahedron, and
tetrahedron are $12, 6$, and $4$, respectively.
\end{proof}

\begin{proposition}
\label{primitive3}
Suppose $G$ is a primitive subgroup of $\mathrm{GL}_3(\C)$.
Then $G/Z(G) \cong H/Z(H)$, where $H$ is one of the following:
\begin{enumerate}
    \item \label{prim3:3A6}
    $H = 3.A_6$, the Valentiner group, with $[H:Z(H)] = 360$;
    \item \label{prim3:Hessian}
    $H$ is the triple cover of the Hessian group, which has $[H:Z(H)] = 216$, or is subgroup of the triple cover with index $3$ or $6$;
    \item \label{prim3:Klein}
    $H = \mathrm{PSL}_2(7)$, the Klein group, with $[H:Z(H)] = 168$;
    \item \label{prim3:A5}
    $H = A_5$, the group of rotational symmetries of the dodecahedron in $\mathrm{SO}(3,\R) \subset \mathrm{GL}_3(\C)$, with $[H:Z(H)] = 60$.
\end{enumerate}
The smallest possible degree of a semi-invariant polynomial is $6$ in 
case \eqref{prim3:3A6} and for the Hessian and its index $3$ subgroup 
in case \eqref{prim3:Hessian}.  The smallest degree is $4$ in case \eqref{prim3:Klein}. 
For case \eqref{prim3:A5},
the smallest degrees of semi-invariants are $2$ and $4$.
In appropriate coordinates, the unique degree $4$ semi-invariant of the Klein
group is the Klein quartic
$$x^3 y + y^3 z + z^3 x.$$
The next smallest degree semi-invariant has degree $6$.

The unique degree $6$ semi-invariant of the Valentiner group is the 
Wiman sextic
$$10 x^3 y^3 + 9 (x^5+y^5)z - 45 x^2 y^2 z^2 - 135 x y z^4 + 27 z^6.$$
The next smallest degree semi-invariant has degree $12$.
The unique degree $6$ semi-invariant of the triple cover of the Hessian group is
$$x^6 + y^6 + z^6 - 10(x^3 y^3 + x^3 z^3 + y^3 z^3).$$
\end{proposition}

\begin{proof}
This classification dates back at least to \cite[Chapter V]{Blichfeldt},
though Blichfeldt's list ``misses" the subgroups of $\mathrm{SL}_3(\C)$ 
which are 
the product of $\mathrm{PSL}_2(7)$ or $A_5$ with the scalar subgroup
$\mu_3$.  Blichfeldt's work also includes explicit matrix generators.

The claims about semi-invariant polynomials of $3.A_6$ and 
$\mathrm{PSL}_2(7)$
are well-known and date back to Wiman \cite{Wiman} and Klein \cite{Klein}, 
respectively. 
The triple cover of the Hessian group, which has order $648$, is the 
extension of the nonabelian group $H_3$ of order $27$ and exponent $3$ by 
$\mathrm{SL}_2(3) \cong 2.A_4$.  Its image in $\mathrm{PGL}_3(\C)$,
the Hessian group, is the automorphism group of the Hesse configuration in
$\mathbb{P}^2_{\C}$ consisting of $9$ points and $12$ lines. The smallest
degree semi-invariant for this group is $6$ \cite{AD}, and the same is true
for the index $3$ primitive subgroup, by \texttt{MAGMA} computation.
\end{proof}

\begin{proposition}
\label{primitive4}
Suppose $G$ is a primitive subgroup of $\mathrm{GL}_4(\C)$ 
satisfying $[G:Z(G)] \geq 648$.
Then $G/Z(G) \cong H/Z(H)$, where $H$ is one of the following:
\begin{enumerate}
    \item \label{prim4:Sp4(3)} 
    $H = \mathrm{Sp}_4(3)$, with $[H:Z(H)] = 25920$;
    \item \label{prim4:H4}
    $H$ is a primitive subgroup of the normalizer $N$ of an extraspecial group
    $H_4$ of order $32$ in $\mathrm{SL}_4(\C)$, where $[N:Z(N)] = 11520$;
    \item \label{prim4:tensor}
    $H$ is an extension by a group of order at most $2$ of a group $A \times B/Z$, where $A,B$ are among \eqref{prim2:A5}, \eqref{prim2:S4}, \eqref{prim2:A4} of \Cref{primitive2} and $Z$ is a central subgroup of order $2$ contained in neither factor; here $[H:Z(H)] \leq 7200$;
    \item \label{prim4:2A7}
    $H = 2.A_7$, with $[H:Z(H)] = 2520$;
    \item \label{prim4:2S6}
    $H = 2.S_6$, with $[H:Z(H)] = 720$.
    \end{enumerate}
The smallest degree for a polynomial semi-invariant is $12$ for 
\eqref{prim4:Sp4(3)} and $8$ for \eqref{prim4:2A7} and \eqref{prim4:2S6}.  The group $N$ in 
\eqref{prim4:H4} fits into an
exact sequence 
$$0 \rightarrow \tilde{H}_4 \rightarrow N \xrightarrow{\rho} S_6 \rightarrow 0,$$
where $\tilde{H}_4$ is generated by $H_4$ and $iI_4$. The smallest degree 
for a polynomial semi-invariant of a subgroup in case \eqref{prim4:H4} is 
always at least $4$, and is at least $8$ unless $H$ is the subgroup
corresponding to the exceptional embedding of $S_5 \subset S_6$ or its
subgroup of index $2$. For these last two groups,
the unique semi-invariant polynomial of degree $4$ is 
$$x_0^4 + x_1^4 + x_2^4 + x_3^4 - 6(x_0^2 x_1^2 + x_0^2 x_2^2 + 
x_0^2 x_3^2 + x_1^2 x_2^2 + x_1^2 x_3^2 + x_2^2 x_3^2).$$
Finally, in case \eqref{prim4:tensor}, every semi-invariant of degree
less than $6$ is a power of the quadric $x_0 x_3 - x_1 x_2$.  This is 
true for degree less than $12$ if additionally $[H:Z(H)] > 1152$.
\end{proposition}

\begin{proof}
The classification of primitive subgroups of $\mathrm{GL}_4(\C)$
is originally due to Blichfeldt \cite[Chapter VII]{Blichfeldt}; all
the groups above are listed (with generators) in his work.  Note that
the summary of the list in \cite[Section 8.5]{Feit} misses a few groups,
including $2.A_5$ and some of the tensor product cases.  Using the
matrix generators, all the claims about semi-invariant polynomials
may be verified in \texttt{MAGMA}.  We note that the extension
of $\tilde{H}_4$ by $S_5$ in case \eqref{prim4:H4}
with semi-invariant of degree $4$ is labeled ``$19^{\circ}$" in
\cite[p. 172]{Blichfeldt}.
\end{proof}

\begin{proposition}
\label{primitive5}
Suppose $G$ is a primitive subgroup of $\mathrm{GL}_5(\C)$ 
satisfying $[G:Z(G)] \geq 9720$.  Then $G$ is generated by scalar 
matrices and the group $\mathrm{PSp}_4(3) \subset \mathrm{GL}_5(\C)$,
which is simple of order $25920$. 
The smallest degree for a semi-invariant polynomial of this group is $4$.
\end{proposition}

\begin{proof}
The classification in dimension $5$ is due to Brauer 
\cite[Theorem 9A]{Brauer} 
(see also \cite[Section 8.5]{Feit}).  We note that the groups labeled
(IV), (V) in Brauer's theorem are subgroups of 
$N = H_5 \rtimes \mathrm{SL}_2(5)$, which has $[N:Z(N)] = 3000 < 9720$.
The claim on semi-invariants is readily verified using \texttt{GAP}.
\end{proof}

\begin{proposition}
\label{primitive6}
Suppose $G$ is a primitive subgroup of $\mathrm{GL}_6(\C)$ 
satisfying $[G:Z(G)] \geq 174960$.  Then $G/Z(G) \cong H/Z(H)$, where
$H$ is one of the following:
\begin{enumerate}
    \item \label{prim6:U4(3)}
    $H = 6_1.\mathrm{PSU}_4(3).2_2$, with $[H:Z(H)] =  6531840$
    or its subgroup of index $2$;
    \item \label{prim6:J2}
    $H = 2.J_2$, a central extension of the sporadic 
    simple Janko group $J_2$, with $[H:Z(H)] = 604800$.
\end{enumerate}
The smallest degree semi-invariant polynomial in case \eqref{prim6:J2} is 
$12$, while in case \eqref{prim6:U4(3)} it is $6$; the unique
semi-invariant polynomial of that degree is
$$\sum_{i = 0}^5 x_i^6 - 10 \sum_{0 \leq i < j \leq 5} x_i^3 x_j^3 - 180 x_0 x_1 x_2 x_3 x_4 x_5.$$
\end{proposition}

\begin{proof}
The classification in dimension $6$ is due to Lindsay \cite{Lindsay1},
where we simply read off the size of each group.  Note that his cases 
(I) and (II), which deal with tensor products of subgroups of 
$\mathrm{GL}_2(\C)$ and $\mathrm{GL}_3(\C)$ and subgroups thereof, have 
center of index at most $60 \cdot 360 = 21600$, which is too small.
The semi-invariants of $2.J_2$ are readily computable with the character
tables in \texttt{GAP}, while a description of explicit matrix 
generators for $6_1.\mathrm{PSU}_4(3).2_2$
appears in \cite[p. 407]{Lindsay2}.  The sextic polynomial shown was
first calculated by Todd \cite{Todd}.
\end{proof}

\begin{proposition}
\label{primitive8}
Suppose $G$ is a primitive subgroup of $\mathrm{GL}_8(\C)$ satisfying
$[G:Z(G)] \geq 88179840$.  Then $G/Z(G) \cong H/Z(H)$, where $H$
is one of the following:
\begin{enumerate}
    \item \label{prim8:O8(2)}
    $H = 2.O_8^+(2).2$, with $[H:Z(H)] = 348364800$, or 
    its subgroup of index $2$;
    \item \label{prim8:H_8}
    $H$ contains a subgroup isomorphic to a central product $T$ of an 
    extraspecial group $H_8$ of order $2^7$ with $\mu_4$,
    and $H/T$ is a subgroup of $\mathrm{Sp}_6(2)$; such an $H$
    has $[H:Z(H)] \leq 185794560$.
\end{enumerate}
Moreover, any such $G$ has semi-invariant polynomials only in even degree.
\end{proposition}

\begin{proof}
The classification in dimension $8$ is due to \cite{HW} when
$7$ does not divide the order of the group and
\cite[Theorem A]{Feitdim8} in the remaining cases.  Inspection of
these lists reveals that only cases (iii) and (ix) of
\cite[Theorem A]{Feitdim8} have large enough order to be considered here;
these are cases \eqref{prim8:H_8} and \eqref{prim8:O8(2)}, respectively.  
Note
that $2.O_8^+(2).2$ is also the Weyl group $W(E_8)$, and a description of
explicit matrix generators for this group appears in
\cite[p. 407]{Lindsay2}.  In case \eqref{prim8:H_8}, 
$$[H:Z(H)] \leq \frac{(2^7 \cdot 4 /2)|\mathrm{Sp}_6(2)|}{2} = 185794560.$$
Both of these cases yield only even degree semi-invariant polynomials
because any subgroup with abelian quotient contains $-I_8$.
\end{proof}

\section{Background and Main Results}
\label{sect:mainresults}

In this section, we'll state a more precise version of 
\Cref{mainthm}, with a full description of the optimal examples in all
the exceptional cases.  We'll also provide some background and motivation
for this problem.

In order to study the automorphism group of the smooth hypersurface
$X \coloneqq \{f = 0\} \subset \mathbb{P}^{n+1}$ of degree $d$,
we will focus on \textit{linear
automorphisms} of $X$, namely those which extend to the ambient
variety $\mathbb{P}^{n+1}$.  We denote by 
$\mathrm{Lin}(X) \subset \mathrm{Aut}(X)$ the subgroup of all
linear automorphisms of $X$.

A well-known result is that $\mathrm{Aut}(X) = \mathrm{Lin}(X)$
whenever $n \geq 1$, $d \geq 3$, and $(n,d)$ is not $(1,3)$ or
$(2,4)$.  For $n \geq 3$, this is a consequence of the 
Grothendieck-Lefschetz theorem \cite[Expos\'{e} XII, 
Corollaire 3.6]{SGA2}, while the $n = 2$ and $n = 1$ cases are
due to Matsumura-Monsky \cite[Theorem 2]{MM} and Chang 
\cite[Theorem 1]{Chang}. The
group $\mathrm{Lin}(X)$ is furthermore finite for $n \geq 1$
and $d \geq 3$ (this is stated for $n \geq 2$ in 
\cite[Theorem 1]{MM}, and the same proof works for plane curves).
We'll state our main result below in terms of $\mathrm{Lin}(X)$;
this will imply \Cref{mainthm} in light of the above comments, but
is also stronger, since it includes the cases $(n,d) = (1,3), (2,4)$.

For most of the paper, we work with subgroups of $\mathrm{GL}_N(\C)$
instead of $\mathrm{PGL}_N(\C)$.  Therefore, we define
$\mathrm{Lin}(f) = \{g \in \mathrm{GL}_{n+2}(\C): g \cdot f = f\}$.

\begin{lemma}
\label{linf}
The natural homomorphism $\mathrm{Lin}(f) \rightarrow \mathrm{Lin}(X)$
induced by $\mathrm{GL}_{n+2}(\C) \rightarrow \mathrm{PGL}_{n+2}(\C)$
is a central extension of $\mathrm{Lin}(X)$ by the group $\mu_d I$
of scalar matrices.
\end{lemma}

\begin{proof}
For any representative $g \in \mathrm{GL}_{n+2}(\C)$ of an element
in $\mathrm{Lin}(X)$, we must have $g \cdot f = a f$ for some
constant $a \in \C$.  Composing $g$ with an appropriate scalar
transformation yields an element in $\mathrm{Lin}(f)$ with the
same image in $\mathrm{PGL}_{n+2}(\C)$, so 
$\mathrm{Lin}(f) \rightarrow \mathrm{Lin}(X)$ is surjective.
Its kernel consists of the scalar matrices that preserve $f$,
which is precisely $\mu_d I$ since $f$ is homogeneous of 
degree $d$.
\end{proof}

We'll see that the example with largest linear automorphism group in most 
dimensions and degrees is the Fermat hypersurface.

\begin{example}
For any positive integers $n$ and $d$, the \textit{Fermat
hypersurface} of dimension $n$ and degree $d$ is the 
hypersurface $X \subset \mathbb{P}^{n+1}$ defined by
the equation 
$$x_0^d + x_1^d + \cdots + x_{n+1}^d = 0.$$
This hypersurface is smooth and it's straightforward
to see that its automorphism group contains the
semidirect product $(\Z/d)^{\oplus (n+1)} \rtimes S_{n+2}$.
Here $(\Z/d)^{\oplus (n+1)}$ is the subgroup of 
automorphisms given by multiplying each coordinate by
a $d$th root of unity (modulo overall scalar transformations)
and $S_{n+2}$ permutes the coordinates.
When $d \geq 3$ we have $\mathrm{Lin}(X) = (\Z/d)^{\oplus (n+1)} \rtimes S_{n+2}$.
(In fact, the same holds more generally over an 
algebraically closed field of any characteristic $p$, so long as 
$p$ does not divide $d$ and $d-1$ is not a power of $p$ \cite{Shioda,Kontogeorgis}.)
In particular, $|\mathrm{Lin}(X)| = (n+2)!d^{n+1}$. 
Similarly, $\mathrm{Lin}(f) = (\Z/d)^{\oplus (n+2)} \rtimes S_{n+2}$,
with order $|\mathrm{Lin}(f)| = (n+2)!d^{n+2}$.
\end{example}

Significant heuristic evidence suggests that the Fermat
hypersurface should achieve the largest automorphism group,
at least generically.

For instance, a result of Howard and Sommese \cite{HS} shows for each
fixed $n$ that there exists a constant $C_n$ depending only on $n$ with the
following property: for every smooth hypersurface $X$ of dimension $n$ and 
degree $d \geq 3$, $|\mathrm{Lin(X)}| \leq C_n d^{n+1}$.
This shows at least that the Fermat hypersurface example has the right
asymptotics in degree, though the result provides no effective 
procedure for computing $C_n$.

Also, a result of Collins \cite{Collins} shows that the symmetric group
$S_{N+1}$, embedded in $\mathrm{GL}_N(\C)$ via the regular representation,
has the largest Jordan constant
of any subgroup of $\mathrm{GL}_N(\C)$ for $N \geq 71$.
The \textit{Jordan constant} $J$ of the 
finite group $G$ is the smallest index $[G:N]$ of a normal abelian subgroup
$N$ of $G$, and roughly measures ``how far $G$ is from being abelian."
Since the symmetric group is asymptotically ``largest" in this sense, it's
plausible that the Fermat hypersurface has largest automorphism group.
Collins' result doesn't directly prove this, however, even for $N \geq 71$,
and it furthermore suggests that exceptional examples of small dimension could
exist.  The first author asked in \cite[Question 3.9]{Esser} whether any
exceptional examples exist in dimension at least $2$.
Our main result confirms that the Fermat hypersurface has largest
automorphism group for most dimensions and degrees and 
enumerates all exceptions.

\begin{theorem}
\label{mainthmfull}
Suppose that $X \subset \mathbb{P}^{n+1}_{\C}$ is a smooth hypersurface
of degree $d$ with $n \geq 1$ and $d \geq 3$.  The Fermat hypersurface,
with linear automorphism group of order $(n+2)!d^{n+1}$, achieves the
largest $|\mathrm{Lin}(X)|$ in dimension $n$ and degree $d$, except
in the cases $(n,d) = (1,4), (1,6), (2,4), \\
(2,6), (2,12), (4,6),$ and $(4,12)$.

Further, any $X$ with $|\mathrm{Lin}(X)| \geq (n+2)!d^{n+1}$
must be isomorphic to the Fermat hypersurface or one of eight
exceptional examples. These examples and the matrix generators for their 
linear automorphism groups are described in Examples \ref{ex:(1,4)} through
\ref{ex:(4,12)}.
\end{theorem}

\begin{remark}
\begin{enumerate}
\item A few cases of \Cref{mainthmfull} were known previously. When $n = 1$, the theorem was proved independently
by works of Pambianco \cite{Pambianco} and Harui \cite{Harui},
building off many previous results in specific degrees. Several recent works partially or completely classify
the automorphism groups that can occur for other specific 
small values of $n$ and $d$.  As a result,
the Fermat hypersurface was already known to have the largest 
automorphism group for the following pairs
$(n,d)$: $(2,3)$ \cite{DD}, $(3,3)$ \cite{WY}, $(3,5)$
\cite{OY}, $(4,3)$ \cite{LZ}, and $(5,3)$ \cite{YangYuZhucubic}. 
All other cases were open, to
the authors' knowledge.
\item The case $(n,d) = (2,4)$ is not included in the earlier statement 
\Cref{mainthm} because $\mathrm{Aut}(X)$ and $\mathrm{Lin}(X)$ can
differ in this case and some quartic K3's have infinite automorphism group.
\item For each of the pairs listed in \Cref{mainthmfull}, there is exactly one hypersurface, up to isomorphism, with $|\mathrm{Lin}(X)| > (n+2)!d^{n+1}$.  In addition, there is an eighth example, Example \ref{ex:(1,6)-2}, with $(n,d) = (1,6)$,
$|\mathrm{Lin}(X)| = (n+2)!d^{n+1} = 216$, but $X$ is not isomorphic to the Fermat sextic, bringing the total number of exceptional examples in the theorem to eight.
\item Several of the examples of dimension at least $2$ appear in the literature and have other interesting properties, though it was not known previously that they had largest automorphism group.  The degree $6$ and $12$
surfaces in Example \ref{ex:(2,6)} and Example \ref{ex:(2,12)}, respectively, have the most known lines of any surfaces of those degrees \cite{BS}, and the degree $6$ example is one of only a few known surfaces of degree at least $5$ of maximal Picard rank \cite{Heijne}.  The equation of the sextic fourfold in Example \ref{ex:(4,6)} is an invariant of a complex reflection group first calculated by Todd \cite{Todd}.
\end{enumerate}

\end{remark}

We list the optimal hypersurfaces in the exceptional cases below,
including matrix generators for $\mathrm{Lin}(f)$.

\begin{example}[$n = 1, d = 4$]
\label{ex:(1,4)}
Let $X \subset \mathbb{P}_{\C}^2$ be the Klein quartic curve defined by
$$f \coloneqq x_0^3 x_1 + x_1^3 x_2 + x_2^3 x_0 = 0.$$
The automorphism group of $X$ is $\mathrm{PSL}_2(7)$.

Explicitly, $\mathrm{Lin}(f)$ is an extension of $\mathrm{PSL}_2(7)$
by $\mu_4$ generated by the following matrices,
$$M_1 = i\begin{pmatrix}\epsilon^4 & 0 & 0 \\ 0 & \epsilon^2 & 0 \\ 0 & 0 & 
\epsilon \end{pmatrix}, M_2 = \begin{pmatrix}0 & 0 & 1 \\ 1 & 0 & 0 \\ 0 & 1 
& 0 \end{pmatrix}, M_3 = \alpha \begin{pmatrix}\epsilon - \epsilon^6 &
\epsilon^2 - \epsilon^5 & \epsilon^4 - \epsilon^3 \\ \epsilon^2 - \epsilon^5 
& \epsilon^4 - \epsilon^3 & \epsilon - \epsilon^6 \\ \epsilon^4 - \epsilon^3
& \epsilon - \epsilon^6 & \epsilon^2 - \epsilon^5 \end{pmatrix},$$
where $\epsilon$ is a primitive $7$th root of unity and $\alpha = \frac{1}{\sqrt{-7}}$.

This group has order $672$, and its quotient 
$\mathrm{Lin}(X) = 
\mathrm{PSL}_2(7) \subset \mathrm{PGL}_3(\C)$ has order $168$.
\end{example}

\begin{example}[$n = 1, d = 6$, Example 1]
\label{ex:(1,6)}
Let $X \subset \mathbb{P}_{\C}^2$ be the Wiman sextic curve defined by 
$$f \coloneqq 10 x_0^3 x_1^3 + 9 (x_0^5+x_1^5)x_2 - 45 x_0^2 x_1^2 x_2^2 - 135 x_0 x_1 x_2^4 + 27 x_2^6 = 0.$$
The automorphism group of $X$ is isomorphic to $A_6$.  The above equation
is the presentation of the Wiman sextic due to Wiman himself \cite{Wiman}.
The matrix generators for the automorphism
group in these coordinates are very messy.
For completeness, we provide a set of matrix generators for 
$\mathrm{Lin}(f)$ in a different coordinate system:
$$M_1 = \begin{pmatrix} 1 & 0 & 0 \\ 0 & -1 & 0 \\ 0 & 0 & 1 \end{pmatrix},
M_2 = \begin{pmatrix} 0 & 1 & 0 \\ 0 & 0 & 1 \\ 1 & 0 & 0 \end{pmatrix},
M_3 = \begin{pmatrix} 1 & 0 & 0 \\ 0 & 0 & \zeta^2 \\ 0 & -\zeta & 0 \end{pmatrix},
M_4 = \frac{1}{2}\begin{pmatrix} 1 & \tau & \tau^{-1} \\ \tau^{-1} & \tau & 1 \\ \tau & -1 & \tau^{-1} \end{pmatrix},$$
where $\zeta$ is a primitive third root of unity and $\tau = \frac{1+\sqrt{5}}{2}.$

The group $\mathrm{Lin}(f)$ has order $2160$, and its quotient 
$\mathrm{Lin}(X) = A_6 \subset \mathrm{PGL}_3(\C)$ has order $360$.
\end{example}

\begin{example}[$n = 1, d = 6$, Example 2]
\label{ex:(1,6)-2}
Let $X \subset \mathbb{P}_{\C}^2$ be the sextic curve defined by 
$$f \coloneqq x_0^6 + x_1^6 + x_2^6 - 10(x_0^3 x_1^3 + x_0^3 x_2^3 + x_1^3 x_2^3) = 0.$$
The automorphism group of $X$ is isomorphic to the Hessian group
of order $216$ in $\mathrm{PGL}_3(\C)$.  Explicitly, $\mathrm{Lin}(f)$ is
the extension of the Hessian group by $\mu_6$ generated by the following matrices:
$$M_1 = \begin{pmatrix} 1 & 0 & 0 \\ 0 & \omega & 0 \\ 0 & 0 & \omega^2 \end{pmatrix}, M_2 = \begin{pmatrix} 0 & 1 & 0 \\ 0 & 0 & 1 \\ 1 & 0 & 0 \end{pmatrix}, M_3 = \frac{1}{i \sqrt{3}}\begin{pmatrix} 1 & 1 & 1 \\ 1 & \omega & \omega^2 \\ 1 & \omega^2 & \omega \end{pmatrix}, M_4 = \begin{pmatrix} \epsilon & 0 & 0 \\ 0 & \epsilon & 0 \\ 0 & 0 & \epsilon^5 \end{pmatrix},$$
where $\omega$ is a primitive third root of unity and $\epsilon$
is a primitive sixth root of unity with $\epsilon^2 = \omega$.

The group $\mathrm{Lin}(f)$ has order $1296$, and its quotient 
$\mathrm{Lin}(X) \subset \mathrm{PGL}_3(\C)$ has order $216$.
\end{example}

\begin{example}[$n = 2, d = 4$]
\label{ex:(2,4)}
Let $X \subset \mathbb{P}_{\C}^3$ be the smooth surface defined by the equation
$$f \coloneqq x_0^4 + x_1^4 + x_2^4 + x_3^4 - 6(x_0^2 x_1^2 + x_0^2 x_2^2 + x_0^2 x_3^2 + x_1^2 x_2^2 + x_1^2 x_3^2 + x_2^2 x_3^2) = 0.$$
The linear automorphism group of $X$ fits into an exact sequence
$$0 \rightarrow \mu_2^{\oplus 4} \rightarrow \mathrm{Lin}(X) \rightarrow S_5 
\rightarrow 0,$$
More explicitly, the group $\mathrm{Lin}(f)$ is generated by the matrices
$$M_1 = \begin{pmatrix}
    0 & 0 & 1 & 0 \\ 0 & 0 & 0 & 1 \\ 1 & 0 & 0 & 0 \\ 0 & 1 & 0 & 0
\end{pmatrix},
M_2 = \begin{pmatrix}
    0 & 1 & 0 & 0 \\ 1 & 0 & 0 & 0 \\ 0 & 0 & 0 & 1 \\ 0 & 0 & 1 & 0
\end{pmatrix},
M_3 = \begin{pmatrix}
    1 & 0 & 0 & 0 \\ 0 & -1 & 0 & 0 \\ 0 & 0 & -1 & 0 \\ 0 & 0 & 0 & 1
\end{pmatrix},$$
$$M_4 = \begin{pmatrix}
    1 & 0 & 0 & 0 \\ 0 & 1 & 0 & 0 \\ 0 & 0 & -1 & 0 \\ 0 & 0 & 0 & -1
\end{pmatrix},
M_5 = \frac{1+i}{2}\begin{pmatrix}
    -i & 0 & 0 & i \\ 0 & 1 & 1 & 0 \\ 1 & 0 & 0 & 1 \\ 0 & -i & i & 0
\end{pmatrix},
M_6 = \begin{pmatrix}
    1 & 0 & 0 & 0 \\ 0 & 1 & 0 & 0 \\ 0 & 0 & 1 & 0 \\ 0 & 0 & 0 & -1
\end{pmatrix}.$$
The subgroup of $\mathrm{GL}_4(\C)$ generated by these matrices has order
$7680$ with scalar subgroup of order $4$.  The image of this group in
$\mathrm{PGL}_4(\C) = \mathrm{Aut}(\mathbb{P}^3)$ has order
$|\mathrm{Lin}(X)| = 1920$.

\end{example}

\begin{example}[$n = 2, d = 6$]
\label{ex:(2,6)}
Let $X \subset \mathbb{P}^3_{\C}$ be
the smooth surface defined by the equation
$$f \coloneqq x_0^5 x_1 - x_0 x_1^5 + x_2^5 x_3 - x_2 x_3^5 = 0.$$
The automorphism group of $X$ is the wreath product of a 
central extension of $S_4$ by $\mu_6$ with $S_2$, modulo its center.

Explicitly, let 
$$M_1 = \frac{1}{2}\begin{pmatrix} 1+i & 1+i \\ -1+i & 1-i \end{pmatrix},
M_2 = \begin{pmatrix} 0 & 1 \\ -1 & 0 \end{pmatrix},
M_3 = \begin{pmatrix} \epsilon & 0 \\ 0 & \epsilon^{19} \end{pmatrix},$$
where $\epsilon$ is a primitive $24$th root of unity with $\epsilon^6 = i$.

The matrices $M_1,M_2,M_3$ generate a central extension of $S_4$ by $\mu_6$ in
$\mathrm{GL}_2(\C)$.  The group $\mathrm{Lin}(f)$ for $f$ above is generated
by $4 \times 4$ block matrices of the form 
$$\begin{pmatrix} M_i & 0 \\ 0 & M_j \end{pmatrix},$$
for $i,j \in \{1,2,3\}$ and the block matrix
$$\begin{pmatrix} 0 & I_2 \\ I_2 & 0 \end{pmatrix}.$$

The subgroup of $\mathrm{GL}_4(\C)$ generated by these matrices
has order $41472$, with  scalar subgroup of order $6$. 
The image of this group in $\mathrm{PGL}_4(\C) = \mathrm{Aut}(\mathbb{P}^3)$
has order $|\mathrm{Lin}(X)| = 6912$.
\end{example}

\begin{example}[$n = 2, d = 12$]
\label{ex:(2,12)}
Let $X \subset \mathbb{P}^3_{\C}$ be the smooth surface defined by the
equation
$$f \coloneqq x_0 x_1 (x_0^{10} + 11 x_0^5 x_1^5 - x_1^{10}) + 
x_2 x_3 (x_2^{10} + 11 x_2^5 x_3^5 - x_3^{10}) = 0.$$
The automorphism group of $X$ is the wreath product of a central
extension of $A_5$ by $\mu_{12}$ with $S_2$, modulo its center.

Explicitly, let
$$M_1 = \begin{pmatrix} \zeta &  0 \\ 0 & \zeta^{49} \end{pmatrix},
M_2 = \frac{1}{\sqrt{5}}\begin{pmatrix} -\epsilon + \epsilon^4 & 
\epsilon^2 - \epsilon^3 \\ \epsilon^2 - \epsilon^3 & 
\epsilon - \epsilon^4, \end{pmatrix},$$
where $\zeta$ is a primitive $60$th root of unity and 
$\epsilon = \zeta^{12}$.

The matrices $M_1,M_2$ generate a central extension of $A_5$ by 
$\mu_{12}$ in $\mathrm{GL}_2(\C)$. The group $\mathrm{Lin}(f)$
for $f$ above is generated
by $4 \times 4$ block matrices of the form 
$$\begin{pmatrix} M_i & 0 \\ 0 & M_j \end{pmatrix},$$
for $i,j \in \{1,2\}$ and the block matrix
$$\begin{pmatrix} 0 & I_2 \\ I_2 & 0 \end{pmatrix}.$$

The subgroup of $\mathrm{GL}_4(\C)$ generated by these matrices
has order $1036800$, with scalar subgroup of order $12$. 
The image of this group in $\mathrm{PGL}_4(\C) = \mathrm{Aut}(\mathbb{P}^3)$
has order $|\mathrm{Lin}(X)| = 86400$.
\end{example}

\begin{example}[$n = 4, d = 6$]
\label{ex:(4,6)}
Let $X \subset \mathbb{P}^5_{\C}$ be the smooth hypersurface defined by the 
equation
$$f \coloneqq \sum_{i = 0}^5 x_i^6 - 10 \sum_{0 \leq i < j \leq 5} x_i^3 x_j^3 - 180 x_0 x_1 x_2 x_3 x_4 x_5 = 0.$$

The automorphism group of $X$ is an extension of
$\mathrm{PSU}_4(3)$ by $\mu_2$.  Explicitly, $\mathrm{Lin}(f)$ is generated by
the following matrices: all $6 \times 6$ permutation matrices, all
$6 \times 6$ matrices of order $3$ and determinant $1$, and the matrix
$I_6 - Q/3$, where $Q$ has all entries equal to $1$.

The subgroup of $\mathrm{GL}_6(\C)$ generated by these matrices has order
$39191040$, with scalar subgroup of order $6$.  The image of this
group in $\mathrm{PGL}_6(\C) = \mathrm{Aut}(\mathbb{P}^5)$ has order
$|\mathrm{Lin}(X)| = 6531840$.
\end{example}

\begin{example}[$n = 4, d = 12$]
\label{ex:(4,12)}
Let $X \subset \mathbb{P}^5_{\C}$ be the smooth hypersurface defined by the
equation
$$x_0 x_1 (x_0^{10} + 11 x_0^5 x_1^5 - x_1^{10}) + 
x_2 x_3 (x_2^{10} + 11 x_2^5 x_3^5 - x_3^{10}) + 
x_4 x_5 (x_4^{10} + 11 x_4^5 x_5^5 - x_5^{10}) = 0.$$
The automorphism group of $X$ is the wreath product of a central
extension of $A_5$ by $\mu_{12}$ with $S_3$, modulo its center.

Explicitly, let
$$M_1 = \begin{pmatrix} \zeta &  0 \\ 0 & \zeta^{49} \end{pmatrix},
M_2 = \frac{1}{\sqrt{5}}\begin{pmatrix} -\epsilon + \epsilon^4 & 
\epsilon^2 - \epsilon^3 \\ \epsilon^2 - \epsilon^3 & 
\epsilon - \epsilon^4, \end{pmatrix},$$
where $\zeta$ is a primitive $60$th root of unity and 
$\epsilon = \zeta^{12}$.

The matrices $M_1,M_2$ generate a central extension of $A_5$ by 
$\mu_{12}$ in $\mathrm{GL}_2(\C)$. The group $\mathrm{Lin}(f)$ is generated
by $6 \times 6$ block matrices of the form 
$$\begin{pmatrix} M_i & 0 & 0 \\ 0 & M_j & 0 \\ 0 & 0 & M_k \end{pmatrix},$$
for $i,j,k \in \{1,2\}$ and the block matrices
$$\begin{pmatrix} 0 & I_2 & 0 \\ I_2 & 0 & 0 \\ 0 & 0 & I_2 \end{pmatrix},
\begin{pmatrix} 0 & I_2 & 0 \\ 0 & 0 & I_2 \\ I_2 & 0 & 0 \end{pmatrix}.$$
The subgroup of $\mathrm{GL}_6(\C)$ generated by these matrices
has order $2239488000$, with  scalar subgroup of order $12$. 
The image of this group in $\mathrm{PGL}_4(\C) = \mathrm{Aut}(\mathbb{P}^5)$
has order $|\mathrm{Lin}(X)| = 186624000$.
\end{example}

\section{Bounding the Order of Automorphism Groups}
\label{sect:bounds}

In this section, we prove some general bounds on the
automorphism groups of hypersurfaces.  We work with the
finite group $\mathrm{Lin}(f)$, the group of 
linear automorphisms preserving the equation $f$ of a 
hypersurface $X$. By \Cref{linf}, this only differs in
order from the (linear) automorphism group of $X$
by a factor of the degree $d$.

\subsection{A General Upper Bound Formula}
\label{subsect:upperbound}

Suppose $X \subset \mathbb{P}^{n+1}$ is a smooth hypersurface
of degree $d$ with defining polynomial $f$.
The idea for bounding $|\mathrm{Lin}(f)|$ is to decompose it as a subnormal
series and bound each of the subquotients individually.
To do this, we first introduce
some notation relating to the group $G \coloneqq \mathrm{Lin}(f)$
that will be used for the next few sections.

We set $N \coloneqq n+2$ to be the dimension of the vector space
$V \cong \C^N$ such that $\mathbb{P}(V) \cong \mathbb{P}^{n+1}$ is 
the ambient projective space for $X$.  Thus we get a representation
$G \subset \mathrm{GL}(V)$ of $G = \mathrm{Lin}(f)$.  We can better
understand this representation by decomposing it as follows:

\begin{lemma}
\label{primdecomp}
For any subgroup $G \subset \mathrm{GL}(V)$, there is a direct sum
decomposition $V = V_1 \oplus \cdots \oplus V_r$ with the following properties: 
\begin{enumerate}
    \item $G$ permutes the subspaces $V_1, \ldots, V_r$;
    \item For each $i = 1,\ldots, r$, the subgroup $H_i \coloneqq 
    \mathrm{Stab}_G(V_i)$ of elements sending $V_i$ to itself acts
    primitively on $V_i$.
\end{enumerate}
\end{lemma}

This lemma is fairly standard.  For completeness, we recall the 
proof, which also appears as \cite[Lemma 1]{Collins}.

\begin{proof}
If $G$ is primitive, the trivial decomposition satisfies the conditions of 
the lemma. Suppose instead that $V$ is an irreducible but not primitive 
$G$-representation.  Then there exists an $H_1$-subrepresentation
$V_1 \subset V$ and a proper subgroup $H_1 \subset G$ such that 
$V = V_1^G$ is induced from 
$V_1$.  Choose $H_1$ a minimal subgroup with this property.  Then $H_1$ acts
primitively on $V_1$ and its images under $G$ form the remainder of the 
decomposition in the lemma.  When $V$ is not irreducible, we can repeat
this procedure for every irreducible subrepresentation.
\end{proof}

We call a direct sum $V = V_1 \oplus \cdots \oplus V_r$ a \textit{primitive
decomposition} for $G$ if it satisfies the conditions of \Cref{primdecomp}.

From now on, fix a primitive decomposition $V = V_1 \oplus \cdots \oplus V_r$
for $G = \mathrm{Lin}(f)$. This decomposition has an associated
\textit{dimension partition} $\pi = (N^{\mu_N},\ldots,1^{\mu_1})$, where
each $\mu_k$ is the number of times a subspace of dimension $k$ appears in 
the decomposition and $r = \mu_N + \cdots + \mu_1$.

We also set $H_i \coloneqq \mathrm{Stab}_G(V_i)$, and 
let 
$$H \coloneqq \bigcap_{i = 1}^r H_i$$
be the kernel of the permutation action of $G$ on the spaces $V_i$, so that
there is an embedding of $G/H$ in the symmetric group $S_r$.  Then we can
write the following subnormal series for $G$:
$$1 \trianglelefteq Z(H) \trianglelefteq H \trianglelefteq G,$$
where $Z(H)$ is the center of $H$.  The strategy for bounding $|G|$
will be to individually bound the orders of the subquotients $G/H$,
$H/Z(H)$, and $Z(H)$.  We also define $\hat{H}_i \subset \mathrm{GL}(V_i)$
to be the image of $H_i$ under the projection.
This is called the 
$i$th \textit{primitive constituent} of $G$.

The key to bounding the first two subquotients is a replacement
theorem due to Collins that replaces $G$ with a better behaved subgroup
$\hat{G}$ with the same primitive decomposition.

\begin{theorem}{{{\cite[Theorem 2]{Collins}}}}
\label{Collinsrepl}
For any subgroup $G \subset \mathrm{GL}(V)$ with primitive decomposition
$V = V_1 \oplus \cdots \oplus V_r$, there exists a group $\hat{G} \subset 
\mathrm{GL}(V)$ with the same decomposition and such that, if we define
$$\hat{H} = \bigcap_{i=1}^r \mathrm{Stab}_{\hat{G}}(V_i),$$
then the following properties hold:
\begin{enumerate}
    \item $H \subset \hat{H}$;
    \item $\hat{H} = \hat{H}_1 \times \cdots \times \hat{H}_r$, where each
    $\hat{H}_i$ acts faithfully and primitively on $V_i$ and trivially on 
    all other summands $V_j, j \neq i$; further, the image $\hat{H}_i \subset 
    \mathrm{GL}(V_i)$ coincides with the image of the projection $H_i 
    \rightarrow \mathrm{GL}(V_i)$;
    \item $\hat{G}$ permutes the subgroups $\{\hat{H}_i\}$ under 
    conjugation with the same permutation action as on the $\{V_i\}$ under
    the bijection $\hat{H}_i \mapsto V_i$;
    \item $\hat{G} = \hat{H} \rtimes \hat{K}$, where $\hat{K} \subset S_r$
    is the subgroup given by the direct product of symmetric groups acting
    on the orbits of the $\{V_i\}$ under $G$;
    \item $[\hat{G}:\hat{H}] \geq [G:H]$, $[\hat{H}:Z(\hat{H})] \geq 
    [H:Z(H)]$, and $Z(\hat{H}) \supset Z(H)$.
\end{enumerate}
\end{theorem}

A few of the properties listed in \Cref{Collinsrepl} are not in
the statement of \cite[Theorem 2]{Collins} but are clearly stated in
its proof.

The last ingredient we need is a method for bounding the center
$Z(H)$.  We first show that each matrix in $Z(H)$ is a
\textit{block scalar matrix} with respect to 
$V = V_1 \oplus \cdots \oplus V_r$, meaning that it is 
contained in the direct sum of scalar subgroups of the $\mathrm{GL}(V_i)$.
Indeed, if $\hat{H}$ is as in \Cref{Collinsrepl},
any matrix in $Z(\hat{H})$ is block scalar since the 
representation of each factor $\hat{H}_i$ in $\hat{H}$ on $V_i$ is
primitive, in particular irreducible, so the conclusion follows from
\Cref{irrep_center}.  Since $Z(\hat{H}) \supset Z(H)$ by \Cref{Collinsrepl},
every matrix in $Z(H)$ is also block scalar with respect to the
primitive decomposition.  Conversely, any block scalar matrix
commutes with every element of $H$ because $H \subset 
\mathrm{GL}(V_1) \oplus \cdots \oplus \mathrm{GL}(V_r)$, so $Z(H)$
equals the block scalar subgroup of $H$.

Now we can bound $|Z(H)|$ using a generalization of a result of the
first author \cite[Theorem 3.3, Step 3]{Esser}, which bounded the 
order of an abelian subgroup of $\mathrm{Lin}(f)$:

\begin{theorem}
\label{absubgrp}
Suppose $X = \{f = 0\} \subset \mathbb{P}^{n+1} = \mathbb{P}(V)$ is
a smooth hypersurface of degree $d \geq 3$.  Let 
$V = V_1 \oplus \cdots \oplus V_r$ be any decomposition of
$V$ into a direct sum of positive-dimensional subspaces.  Then the 
intersection of $\mathrm{Lin}(f)$ with the group of 
block scalar matrices with respect to this decomposition
has order at most $d^r$.
\end{theorem}

Note that the decomposition here is arbitrary, and doesn't actually
need to be a primitive decomposition for $\mathrm{Lin}(f)$.  Since
any abelian subgroup is diagonalizable (hence block scalar for some
decomposition into one-dimensional subspaces), 
\cite[Theorem 3.3, Step 3]{Esser} is a special case of \Cref{absubgrp}
for $r = N = n+2$.  The fact from
\Cref{linf} that the intersection of $\mathrm{Lin}(f)$
with the scalar matrices is $\mu_d$ is the special case $r = 1$.

\begin{proof}
The only consequence of the smoothness of $X$ that we need for the
proof is the existence of certain monomials with nonzero coefficients
in $f$.  In particular, for each variable $x_j$, where $j = 0,\ldots,n+1$,
$f$ contains a monomial of the form $x_j^d$ or $x_j^{d-1} x_k$, $k \neq j$,
with nonzero coefficient in $f$.  Otherwise, $X$ would be singular at the 
coordinate point of the variable $x_j$ \cite[Proposition 1.2]{Esser}.

Suppose first that each subspace of the decomposition is $1$-dimensional, so that 
we're bounding the size of the diagonal subgroup $A$ of $\mathrm{Lin}(f)$.
This case was proven in \cite[Theorem 3.3, Step 3]{Esser}, but we review
the proof for the reader's convenience.
Each diagonal automorphism $x_j \mapsto c_j x_j$ in $A$ must
preserve each monomial of $f$ individually, so if 
$$f = \sum_{i=1}^s K_i \prod_{j = 0}^{n+1}x_{j}^{m_{ij}}$$
is invariant, then for each $i$,
\begin{equation}
\label{diagonalaut}
  \prod_{j = 0}^{n+1} c_j^{m_{ij}} = 1.  
\end{equation}

Assemble the exponents $m_{ij}$ into a $s \times (n+2)$ matrix $M$.
The equations above for each $i$ imply that
$(\log |c_0|,\ldots,\log |c_{n+1}|)^{\mathsf{T}} \in \ker(M)$.
We claim that the kernel of $M$ is trivial.  To see this, take the
$(n+2) \times (n+2)$ square minor $B$ of $M$ which corresponds to the 
monomials of the special form $x_j^d$ or $x_j^{d-1} x_k$, where
the monomial for $x_j$ goes in the $j$th row.  Then $B$ is invertible and
in particular has positive determinant bounded above by $d^{n+2}$ 
\cite[Lemma 3.6]{Esser}.  It follows that $\ker(M) = 0$ and each
$c_j$ is a unit complex number $c_j = e^{2 \pi i \theta_j}$.
The condition \eqref{diagonalaut} then gives that 
$M (\theta_0,\ldots,\theta_{n+1})^{\mathsf{T}} \in \Z^s$.
The weaker condition $B (\theta_0,\ldots,\theta_{n+1})^{\mathsf{T}} \in \Z^{n+2}$
is already enough to show that $(\theta_0,\ldots,\theta_{n+1})$ is in the
superlattice of $\Z^{n+2}$ spanned by $B^{-1} e_0,\ldots, B^{-1}e_{n+1}$,
where $e_0,\ldots,e_{n+1}$ are the standard basis vectors of $\Z^{n+2}$.
Since each $\theta_j$ is only defined modulo $\Z$, we have an injection
$A \subset B^{-1} \Z^{n+2}/\Z^{n+2}$.  This quotient has order $\det(B) \leq d^{n+2}$,
completing the argument.

Now suppose $V = V_1 \oplus \cdots \oplus V_r$ is an arbitrary decomposition
and let $A$ be the intersection of $\mathrm{Lin}(f)$ with the block
scalar matrices of this decomposition.  We can define a
linear transformation $V \rightarrow \C^r$ by ``collapsing" the coordinates
in each block.  Indeed, if $x_0,\ldots,x_{\dim(V_1)-1}$ are the coordinates
of the first block, map $x_0 \mapsto \beta_0 y_1, 
\ldots,  x_{\dim(V_1)-1} \mapsto
\beta_{\dim(V_1)-1} y_1$, where $\beta_0, \ldots,\beta_{\dim(V_1)-1}$ are 
nonzero constants. After doing similarly on the 
other blocks, $V \rightarrow \C^r$ is
$A$-equivariant for a diagonal action of $A$ on $\C^r$.  Furthermore,
the transformed polynomial $\bar{f}(y_1,\ldots,y_m)$ is $A$-invariant.
It's not necessarily true that $\bar{f}$ defines a smooth hypersurface
any longer, but for a general
choice of constants $\beta_j$, $\bar{f}$ will 
still contain monomials of the form
$y_j^d$ or $y_j^{d-1} y_k$ with nonzero
coefficients for each $j = 1,\ldots,r$.
This is because monomials of this type collapse
to monomials of this type, and 
different monomials will not cancel for general $\beta_j$.
Now, using the same proof as above, we have $|A| \leq d^r$, as required.
\end{proof}

We are now ready to prove a general bound on the order of $\mathrm{Lin}(f)$:

\begin{theorem}
\label{boundlinf}
Using the notation from the beginning of \Cref{subsect:upperbound}:
\begin{equation}
\label{boundlinfeq}
    |\mathrm{Lin}(f)| \leq \mu_1 ! \mu_2 ! \cdots \mu_N ! \prod_{i = 1}^r [\hat{H}_i:Z(\hat{H}_i)]
d^r \leq \mu_1 ! \mu_2 ! \cdots \mu_N ! \prod_{i = 1}^r \Xi(\dim(V_i)) d^r.
\end{equation}

More precisely, 
\begin{align*}
    & |Z(H)| \leq d^r, \\
    & [H:Z(H)] \leq \prod_{i = 1}^r [\hat{H}_i:Z(\hat{H}_i)] \leq \prod_{i = 1}^r \Xi(\dim(V_i)), \text{and } \\
    & [G:H] \leq  \mu_1 ! \mu_2 ! \cdots \mu_N !.
\end{align*}
\end{theorem}

Recall that the $\Xi$ function has values which are the upper bounds for
$[G:Z(G)]$ over primitive linear groups $G$ in a particular dimension (see 
\Cref{Xifunction}).

\begin{proof}
The upper bound for $\mathrm{Lin}(f)$ will be the product of individual
upper bounds for the three numbers $[G:H], [H:Z(H)]$, and $|Z(H)|$.
First, since $Z(H)$ is the block scalar subgroup of $\mathrm{Lin}(f)$, 
it follows from \Cref{absubgrp} that $|Z(H)| \leq d^r$.

Next, $[H:Z(H)] \leq [\hat{H}:Z(\hat{H})]$ by \Cref{Collinsrepl}, and
$[\hat{H}:Z(\hat{H})] = \prod_{i = 1}^r [\hat{H}_i:Z(\hat{H}_i)]$,
which is the expression that appears on the right of the first inequality.
Since each $\hat{H}_i$ acts primitively and faithfully on $V_i$, 
the quantity $[\hat{H}_i:Z(\hat{H}_i)]$ is bounded above by $\Xi(\dim(V_i))$
by definition.

Finally, we know from \Cref{Collinsrepl} that $[G:H] \leq [\hat{G}:\hat{H}]$,
and $\hat{K} = \hat{G}/\hat{H}$ is the subgroup of $S_r$ given by the
direct product of symmetric groups of orbits of the $\{V_i\}$ under $G$.
Of course, the order of this group is bounded above by $r!$, but we can 
do better: it is in fact a subgroup of $S_{\mu_N} \times \cdots \times 
S_{\mu_1} \subset S_r$, where each factor acts only on the subset of 
$\{V_i\}$ of subspaces of the given dimension.  This is because a 
subspace of dimension $k$ can only be permuted by $G$ to a subspace of a
same dimension.  Hence $[G:H] \leq \mu_1 ! \cdots \mu_N!$, completing the
proof.
\end{proof}

\subsection{Exceptional Partitions}
\label{subsec:expart}

Notice that the rightmost expression in \eqref{boundlinfeq} depends
only on the partition $\pi$ of $N$ and the degree $d$.  We define
$$B(\pi,d) \coloneqq \mu_1 ! \mu_2 ! \cdots \mu_N ! \prod_{i = 1}^r \Xi(\dim(V_i)) d^r.$$
An important special case of the above is $\pi = (1^N)$.  When 
$X = \{f = 0\}$ is a Fermat hypersurface, $G = \mathrm{Lin}(f)$
has this $\pi$ as the associated partition
and $G$ realizes the equality $|G| = N! d^N = B((1^N),d)$.

For two partitions $\pi$ and $\pi'$ of the integers $N$ and $N'$, 
respectively, the \textit{concatenation} $\pi \oplus \pi'$ of $\pi$ and 
$\pi'$ is the partition of $N + N'$ whose set of blocks is given by the 
disjoint union of sets of blocks of $\pi$ and $\pi'$.
The following lemma summarizes a few simple properties of the function 
$B(\pi,d)$:

\begin{lemma}
\label{Bfunction}
Let $\pi$ be a partition of $N$ and $\pi'$ a partition of $N'$.  Then:
\begin{enumerate}
    \item \label{Bfunction:1} $B(\pi,d) B(\pi',d) \leq B(\pi \oplus \pi',d)$, and
    \item \label{Bfunction:2} If the sets of block sizes of $\pi$ and $\pi'$ are disjoint,
    $B(\pi,d) B(\pi',d) = B(\pi \oplus \pi',d)$.
\end{enumerate}
\end{lemma}

\begin{proof}
The term $d^r \prod_{i = 1}^r \Xi(\dim(V_i))$ of $B(\pi,d)$ is 
clearly multiplicative for partition concatenation, so the 
inequality in \eqref{Bfunction:1} follows from noticing that $a! b! \leq (a+b)!$ for
any nonnegative integers $a$ and $b$.  In the case that $\pi$ has 
no block of the same size as any block of $\pi'$, the number of
blocks of any given size is $0$ for either $\pi$ or $\pi'$ and
the inequality of factorials becomes an equality.
\end{proof}

Using this bound, we can already prove that the Fermat hypersurface always
has the largest automorphism group of any hypersurface if the
dimension is large enough:

\begin{theorem}
\label{highdim}
Fix any positive integer $N \geq 27$ and degree $d \geq 3$.  Then
$B(\pi,d) < B((1^N),d)$ for any partition $\pi$ of $N$ such that 
$\pi \neq (1^N)$.
\end{theorem}

When $X = \{f = 0\} \subset \mathbb{P}^{n+1}$ is the Fermat hypersurface, 
$\mathrm{Lin}(f) = B((1^N),d)$, so \Cref{highdim} and \Cref{boundlinf} 
together already show that 
the Fermat hypersurface achieves the maximum order automorphism group
in every degree $d \geq 3$ and dimension $n = N-2 \geq 25$.

\begin{proof}
Recall the definition
$$B(\pi,d) = \mu_1 ! \mu_2 ! \cdots \mu_N !
\prod_{i = 1}^r \Xi(\dim(V_i)).$$

The strategy is to show that
replacing blocks of $\pi$ with blocks of dimension $1$ always
increases the quantity $B(\pi,d)$
apart from a small collection of cases.  In particular,
we'll show that these unusual cases only occur when $N < 27$.

For a given partition $\pi$ of $N$, we can write
$\pi = \pi_1 \oplus \cdots \oplus \pi_N$, where the partition $\pi_i$ contains
only the blocks of size $i$.  By \Cref{Bfunction}\eqref{Bfunction:2}, 
$B(\pi,d) = B(\pi_1,d) \cdots B(\pi_N,d)$.  We will successively
replace the partitions $\pi_i$ with partitions 
with blocks of smaller sizes in such a way that the modified
$B(\pi,d)$ is larger.  Indeed, if $\tau$ is a refinement of $\pi_i$, 
$\pi'$ is the new partition
of $N$ with $\pi_i$ replaced by $\tau$, and $B(\pi_i,d) \leq B(\tau,d)$,
then
$$B(\pi,d) = B(\pi_1,d) \cdots B(\pi_i,d) \cdots B(\pi_N,d) \leq 
B(\pi,d) \cdots B(\tau,d) \cdots B(\pi_N,d) \leq B(\pi',d).$$
This argument requires justifying many inequalities involving the
function $B(\pi,d)$.  Note that in each modification, we replace
\textit{all} the blocks of a given size with blocks of strictly smaller size.
In each of these inequalities, there are
more blocks in the new partition than the old, and thus $B(\pi,d)$
has a higher power of $d$ on the right.  This allows us to simplify 
matters and assume that $d$ takes the smallest value $d = 3$.
Also, nearly all the inequalities contain a parameter $k$ or $\ell$.  In
each one, it's straightforward to justify the inequality by
verifying the base case where the parameter takes the smallest 
indicated value and then using induction.

We begin by noting that 
$k!(i+1)! 3^k  < (ki)! 3^{ki} = B((1^{ki}),3)$
holds for all $k \geq 1$, $i > 2$. The left-hand side of this 
inequality equals $B((i^k),3)$ unless $i$ is one of the 
special values from \Cref{primgroups}.
Therefore, if $\pi$ contains
$k$ blocks of size $i \neq 2,3,4,5,6,7,8,9$, or $12$, we modify $\pi$ 
by replacing each block of size $i$ by $i$ blocks of size $1$. By the above,
this modification increases $B(\pi,d)$.  For $i = 12,9,7$,
we claim that performing the same modification increases $B(\pi,d)$.  This
is justified by the following inequalities, which hold for all positive
integers $k$:
\begin{align*}
   & B((12^k),3) =  k! 448345497600^k 3^k < (12k)! 3^{12k} = B((1^{12k}),3);
   \\
   & B((9^k),3) = k! 4199040^k 3^k < (9k)! 3^{9k} = B((1^{9k}),3); \\
   & B((7^k),3) = k! 1451520^k 3^k < (7k)! 3^{7k} = B((1^{7k}),3). \\
\end{align*}
Next, for $i = 8,6,5$, we claim that it's always advantageous 
to replace each block of size $i$ with blocks of size $2$, 
plus one block of size $1$ if $i$ is odd.  This is justified by the following inequalities,
which hold for all positive integers $k$:
\begin{align*}
   & B((8^k),3) = k! 348364800^k 3^k < (4k)! 60^{4k} 3^{4k} = B((2^{4k}),3); \\
   & B((6^k),3) = k! 6531840^k 3^k < (3k)! 60^{3k} 3^{3k}  = B((2^{3k}),3); \\
   & B((5^k),3) = k! 25920^k 3^k < k! (2k)! 60^{2k} 3^{3k} = B((2^{2k},1^k),3). 
\end{align*}
This leaves us with a partition $\pi$ containing only blocks
of sizes $1,2,3$, and $4$.  We can perform various additional 
reductions as follows.  By the inequalities below for the indicated values of $k$ or $\ell$, 
we can replace the blocks of sizes $4$ and $3$ with blocks 
of sizes $2$ and $1$ unless there is at most one block of size $4$ or $3$.
\begin{align*}
    & B((4^k),3) = k! 25920^k 3^k < (2k)! 60^{2k} 3^{2k} = B((2^{2k}),3), k \geq 2; \\
    & B((3^{2 \ell}),3) = (2 \ell)! 360^{2 \ell} 3^{2 \ell}  < (3 \ell)! 60^{3 \ell} 3^{3 \ell} = B((2^{3 \ell}),3), \ell \geq 1; \\
    & B((3^{2 \ell + 1}),3) = (2 \ell + 1)! 360^{2 \ell + 1} 3^{2 \ell + 1} < (3 \ell + 1)! 60^{3 \ell + 1} 3^{3 \ell + 2} = B((2^{3 \ell + 1},1),3), \ell \geq 1.
\end{align*}
At this point $\pi$ only has blocks of sizes $4,3,2,1$, with at most one
each of sizes $4$ and $3$. 
If $N \geq 8$, we can eliminate any lone remaining $4$ and $3$ blocks,
due to the following inequalities, and the fact that $\pi$ would have
to contain a subpartition of one of the following types: $(4,3), (4,2^k),
k \geq 1, (4,1^k), k \geq 4, (3,2^k), k \geq 2, (3,1^k), k \geq 2$.  
Moreover, we can arrange so that this subpartition contains all blocks 
of $\pi$ of the indicated sizes, so that \Cref{Bfunction}\eqref{Bfunction:2} applies.
\begin{align*}
    & B((4,3),3) = 25920 \cdot 360 \cdot 3^2 < 3! 60^3 3^4 = B((2^3,1),3); \\
    & B((4,2^k),3) = k! 25920 \cdot 60^k 3^{k+1} < (k+2)! 60^{k+2} 3^{k+2} = B((2^{k+2}),3), k \geq 1; \\
    & B((4,1^k),3) = k! 25920 \cdot 3^{k+1} < (k+4)! 3^{k+4} = B((1^{k+4}),3), k \geq 4; \\
    & B((3,2^k),3) = k! 360 \cdot 60^k 3^{k+1} < (k+1)! 60^{k+1} 3^{k+2} = B((2^{k+1},1),3), k \geq 2; \\ 
    & B((3,1^k),3) = k! 360 \cdot 3^{k+1} < (k+3)! 3^{k+3} = B((1^{k+3}),3), k \geq 2.
\end{align*}
Finally, under the assumption that $N \geq 8$, we've now reached the point
that any $\pi$ has been replaced by a new $\pi$ with a larger $B(\pi,d)$
such that $\pi$ contains only $1$'s and $2$'s.  The last step is to show
that if $N$ is large enough, we can replace all the $2$'s with $1$'s.
Note that if the number of $2$'s is at least $14$, we can replace all
$2$'s with $1$'s:
$$B((2^k),3) = k! 60^k 3^k < (2k)! 3^{2k} = B((1^{2k}),3), k \geq 14.$$
This means we can arrange for $\pi$ to have at most $13$ $2$'s.  For each
$1 \leq k \leq 13$, one can verify:
$$B((2^k,1^{\ell}),3) = \ell ! k ! 60^k 3^{k + \ell} <
(2k + \ell)! 3^{2k + \ell} = B((1^{2k + \ell}),3), \ell \geq 27-2k.$$
For this last inequality, one takes $\ell = 27-2k$ as the base case
and argues by induction on $\ell$ for each fixed $k$.  In summary, we've shown
that for $N \geq 27$ and any $d \geq 3$, $B((1^N),d)$ is strictly larger than
$B(\pi,d)$, where $\pi$ is any other partition of $N$.
\end{proof}

The hypersurface $X = \{f = 0\}$ of degree $d \geq 3$ can
only have an automorphism group of order equal to or 
larger than that of the Fermat
hypersurface of the same dimension and degree if $G \coloneqq \mathrm{Lin}(f)$
has an associated partition $\pi$ with $B(\pi,d) \geq B((1^{N}),d)$.
We say that such a $\pi \neq (1^N)$ is \textit{exceptional for degree $d$}.
In light of \Cref{highdim}, $\pi$ must be a partition of $N \leq 26$.
A computer search of all such partitions yields exactly $80$ with 
the property that $B(\pi,3) \geq B((1^{N}),3)$ and 
$\pi \neq (1^{N})$; we simply call these
\textit{exceptional partitions}.  They represent the remaining cases
that need to be checked for large automorphism groups.
\Cref{exceptionalpart} lists these partitions, along
with the maximum $d$ for which each is exceptional, and the 
ratio $B/B_F \coloneqq B(\pi,3)/B((1^N),3)$.

\begin{table}
\begin{tabular}[t]{|c|c|c|c|c|}
Number & $N$ & Partition & Max $d$ & $B/B_F$ \\
\hline
No.  1  &  2  &  $(2)$ & 30 & 10.0  \\
\hline
No.  2  &  3  &  $(3)$ & 7 & 6.66  \\
\hline
No.  3  &  3  &  $(2,1)$ & 10 & 3.33  \\
\hline
No.  4  &  4  &  $(4)$ & 10 & 40.0  \\
\hline
No.  5  &  4  &  $(3,1)$ & 3 & 1.67  \\
\hline
No.  6  &  4  &  $(2^2)$ & 17 & 33.3  \\
\hline
No.  7  &  4  &  $(2,1^2)$ & 5 & 1.67  \\
\hline
No.  8  &  5  &  $(5)$ & 3 & 2.67  \\
\hline
No.  9  &  5  &  $(4,1)$ & 6 & 8.00  \\
\hline
No.  10  &  5  &  $(3,2)$ & 5 & 6.67  \\
\hline
No.  11  &  5  &  $(2^2,1)$ & 7 & 6.67  \\
\hline
No.  12  &  5  &  $(2,1^3)$ & 3 & 1.00 \\
\hline
No.  13  &  6  &  $(6)$ & 6 & 37.3  \\
\hline
No.  14  &  6  &  $(4,2)$ & 6 & 26.7  \\
\hline
No.  15  &  6  &  $(4,1^2)$ & 4 & 2.67  \\
\hline
No.  16  &  6  &  $(3^2)$ & 4 & 4.45  \\
\hline
No.  17  &  6  &  $(3,2,1)$ & 3 & 1.11  \\
\hline
No.  18  &  6  &  $(2^3)$ & 12 & 66.7  \\
\hline
No.  19  &  6  &  $(2^2,1^2)$ & 4 & 2.22  \\
\hline
No.  20  &  7  &  $(6,1)$ & 4 & 5.33  \\
\hline
No.  21  &  7  &  $(5,2)$ & 3 & 1.27  \\
\hline
No.  22  &  7  &  $(4,3)$ & 4 & 7.62  \\
\hline
No.  23  &  7  &  $(4,2,1)$ & 4 & 3.81  \\
\hline
No.  24  &  7  &  $(4,1^3)$ & 3 & 1.14  \\
\hline
No.  25  &  7  &  $(3,2^2)$ & 4 & 6.35  \\
\hline
No.  26  &  7  &  $(2^3,1)$ & 6 & 9.53  \\
\hline
No.  27  &  8  &  $(8)$ & 3 & 3.95  \\
\hline
No.  28  &  8  &  $(6,2)$ & 4 & 13.3  \\
\hline
No.  29  &  8  &  $(6,1^2)$ & 3 & 1.33  \\
\hline
No.  30  &  8  &  $(4^2)$ & 5 & 45.7  \\
\hline
No.  31  &  8  &  $(4,2^2)$ & 5 & 19.0  \\
\hline
No.  32  &  8  &  $(3^2,2)$ & 3 & 1.59  \\
\hline
No.  33  &  8  &  $(2^4)$ & 9 & 95.2  \\
\hline
No.  34  &  8  &  $(2^3,1^2)$ & 4 &  2.38  \\
\hline
No.  35  &  9  &  $(6,3)$ & 3 & 2.96  \\
\hline
No.  36  &  9  &  $(6,2,1)$ & 3 &  1.48  \\
\hline
No.  37  &  9  &  $(4^2,1)$ & 3 & 5.08  \\
\hline
No.  38  &  9  &  $(4,3,2)$ & 3 & 2.12  \\
\hline
No.  39  &  9  &  $(4,2^2,1)$ & 3 & 2.12  \\
\hline
No.  40  &  9  &  $(3^3)$ & 3 & 1.06  \\
\hline \end{tabular} \qquad \begin{tabular}[t]{|c|c|c|c|c|}
Number & $N$ & Partition & Max $d$ & $B/B_{\on{F}}$ \\
\hline
No.  41  &  9  &  $(3,2^3)$ & 4 & 5.29  \\ 
\hline
No.  42  &  9  &  $(2^4,1)$ & 5 & 10.6  \\
\hline
No.  43  &  10  &  $(6,4)$ & 3 & 7.11  \\
\hline
No.  44  &  10  &  $(6,2^2)$ & 3 & 5.93  \\
\hline
No.  45  &  10  &  $(4^2,2)$ & 4 & 10.2  \\
\hline
No.  46  &  10  &  $(4^2,1^2)$ & 3 & 1.02  \\
\hline
No.  47  &  10  &  $(4,2^3)$ & 4 &  12.7  \\
\hline
No.  48  &  10  &  $(2^5)$ & 7 &  106  \\
\hline
No.  49  &  10  &  $(2^4,1^2)$ & 3 &  2.12  \\
\hline
No.  50  &  11  &  $(4^2,3)$ & 3 & 1.85  \\
\hline
No.  51  &  11  &  $(4,2^3,1)$ & 3 & 1.15  \\
\hline
No.  52  &  11  &  $(3,2^4)$ & 3 & 3.85  \\
\hline
No.  53  &  11  &  $(2^5,1)$ & 4 & 9.62  \\
\hline
No.  54  &  12  &  $(6^2)$ & 3 & 3.02  \\
\hline
No.  55  &  12  &  $(6,4,2)$ & 3 & 1.08  \\
\hline
No.  56  &  12  &  $(6,2^3)$ & 3 & 2.70  \\
\hline
No.  57  &  12  &  $(4^3)$ & 3 & 11.1  \\
\hline
No.  58  &  12  &  $(4^2,2^2)$ & 3 & 3.08  \\
\hline
No.  59  &  12  &  $(4,2^4)$ & 4 & 7.70  \\
\hline
No.  60  &  12  &  $(2^6)$ & 6 & 96.2  \\
\hline
No.  61  &  12  &  $(2^5,1^2)$ & 3 & 1.60  \\
\hline
No.  62  &  13  &  $(3,2^5)$ & 3 & 2.47  \\
\hline
No.  63  &  13  &  $(2^6,1)$ & 4 & 7.40  \\
\hline
No.  64  &  14  &  $(6,2^4)$ & 3 &  1.18  \\
\hline
No.  65  &  14  &  $(4^3,2)$ & 3 & 1.22  \\
\hline
No.  66  &  14  &  $(4^2,2^3)$ & 3 & 1.01  \\
\hline
No.  67  &  14  &  $(4,2^5)$ & 3 & 4.23  \\
\hline
No.  68  &  14  &  $(2^7)$ & 5 & 74.0  \\
\hline
No.  69  &  14  &  $(2^6,1^2)$ & 3 & 1.06  \\
\hline
No.  70  &  15  &  $(3,2^6)$ & 3 & 1.41  \\
\hline
No.  71  &  15  &  $(2^7,1)$ & 3 & 4.93  \\
\hline
No.  72  &  16  &  $(4,2^6)$ & 3 & 2.11  \\
\hline
No.  73  &  16  &  $(2^8)$ & 4 & 49.3  \\
\hline
No.  74  &  17  &  $(2^8,1)$ & 3 & 2.90  \\
\hline
No.  75  &  18  &  $(2^9)$ & 4 &  29.0  \\
\hline
No.  76  &  19  &  $(2^9,1)$ & 3 & 1.53  \\
\hline
No.  77  &  20  &  $(2^{10})$ & 3 & 15.3  \\
\hline
No.  78  &  22  &  $(2^{11})$ & 3 &  7.27  \\
\hline
No.  79  &  24  &  $(2^{12})$ & 3 & 3.16  \\
\hline
No.  80  &  26  &  $(2^{13})$ & 3 & 1.27 \\ 
\hline
    \end{tabular}
    \hspace{0.2em} 
    \caption{The $80$ exceptional partitions of a positive integer $N$ for which $B(\pi,d) \geq B((1^{N}),d)$ for some $d \geq 3$. The fourth column is the maximum integer $d$ for which $B(\pi,d) \geq B((1^{N}),d)$, which exists since $r < N$. The fifth shows the value $B/B_F \coloneqq B(\pi,3)/B((1^{N}),3)$ to three significant figures.}
    \label{exceptionalpart}
\end{table}

\section{Finishing the Proof of Main Theorem}
\label{sect:proof}

From the results of the last section, it only remains to analyze
the finite number of cases listed in \Cref{exceptionalpart}, all
of which correspond to dimensions $n = N-2$ at most $24$.

Though there are $80$ exceptional partitions, a vast majority of these
will not be realized by hypersurfaces with automorphism group
larger than or equal to that of the Fermat.  This is because the analysis in
\Cref{subsec:expart} does not use any information about invariant theory,
which further constrains the degrees $d$ that can occur for various
partitions.  In this section, we'll finish the proof by 
considering each of the exceptional partitions.

One subtlety is that the exceptional partitions each have some ``wiggle room":
a hypersurface $X = \{f = 0\}$ with $|\mathrm{Lin}(f)|$ larger than or equal to
the order $B((1^N),d) = N! d^N$ of the linear automorphism group of the Fermat
equation does not need to satisfy $|\mathrm{Lin}(f)| = B(\pi,d)$, where
$\pi$ is a partition for a primitive decomposition for $\mathrm{Lin}(f)$.
Indeed, we could instead have $B((1^N),d) \leq |\mathrm{Lin}(f)| < B(\pi,d)$.
To rule out a partition $\pi$,
we must therefore show that $|\mathrm{Lin}(f)|$ always falls short of $B(\pi,d)$
by a factor of more than $B(\pi,d)/B((1^N),d) \leq B(\pi,3)/B((1^N),3)$.
The latter value is the number shown in the fifth column of \Cref{exceptionalpart}.

We'll continue to use the notation from the beginning of 
\Cref{subsect:upperbound}.  So we choose a primitive decomposition
$V = V_1 \oplus \cdots \oplus V_r$ for $G \coloneqq \mathrm{Lin}(f)$ as
in \Cref{primdecomp}, set $H_i \coloneqq \mathrm{Stab}_G(V_i)$, and define
$\hat{H}_i$ to be the $i$th primitive
constituent of $G$, i.e., the image of $H_i$ in $\mathrm{GL}(V_i)$. 
As above, we also let
$$H \coloneqq \bigcap_{i=1}^r H_i.$$
Note that we are using the original group $G$ here rather than its 
replacement $\hat{G}$ from \Cref{Collinsrepl}; replacing the group
can change its invariant polynomials.  Nevertheless, by \Cref{Collinsrepl}(2),
the image of $H_i$ in $\mathrm{GL}(V_i)$ coincides with the $\hat{H}_i$ that
appears in the replacement group. Finally, we'll introduce the notation
$\bar{H}_i$ for the image of $\hat{H}_i$ in $\mathrm{PGL}(V_i)$.
Since $Z(H)$ coincides with the block scalar subgroup of $H$, we have an 
induced inclusion
\begin{equation}
\label{Hinclusion}
H/Z(H) \xhookrightarrow{\phi} \bar{H}_1 \times \cdots \times \bar{H}_r
\subset \mathrm{PGL}(V_1) \oplus \cdots \oplus \mathrm{PGL}(V_r).
\end{equation}
induced by $H \hookrightarrow \hat{H}_1 \times \cdots \times \hat{H}_r$.

The idea of this section is as follows: the upper bound $B(\pi,d)$ for 
$|G|$ is achieved only when each of the three subquotients of $G$ from
\Cref{boundlinf} attains the upper bound for its order.
Conversely, if the order of any
subquotient falls short by some amount, $|G|$ will fall short of 
$B(\pi,d)$ by at least the
same amount.  In particular, if the bound is achieved,
the inclusion in \eqref{Hinclusion} must be an isomorphism.
We use the following argument throughout this section, so we list it here for 
reference:

\begin{setup}
\label{elim_setup}
For some $i$ among $1,\ldots,r$, consider the composition
$$H/Z(H) \xhookrightarrow{\phi} \bar{H}_1 \times \cdots \times \bar{H}_r 
\rightarrow \bar{H}_1 \times \cdots \times \widehat{\bar{H}_i} \times \cdots \times \bar{H}_r.$$
(The second map is the projection to the product omitting $\bar{H}_i$.)
The kernel $\bar{K}$ of this composition can be identified with a subgroup of
$\bar{H}_i$.  Denote by $K$ the preimage of $\bar{K}$ in $H$. 
As a subgroup of
$\mathrm{GL}(V_1) \oplus \cdots \oplus \mathrm{GL}(V_r)$, $K$ maps to
a block scalar subgroup under the projection forgetting $\mathrm{GL}(V_i)$,
because the image of each element of $K$ 
in the projective general linear groups of
the other factors is trivial.  Hence the kernel $K'$ of this projection map
is a subgroup of $H \subset G$ such that (1) $K'$ is a normal subgroup of
$K$ with $K/K'$ abelian; (2) $K'$ acts trivially on $V_j$, where $j \neq i$;
and (3) the image $\bar{K}'$ of $K' \rightarrow \mathrm{PGL}(V_i)$ is a
normal subgroup of $\bar{K}$ with abelian quotient.

Since $f$ is an invariant polynomial of $G$, it is an invariant polynomial
of $K'$.  The group $K'$ acts only on $V_i$ and trivially on other summands,
so setting the variables outside $V_i$ to be arbitrary constants, the 
resulting polynomial $f$ in the variables of $V_i$ must always be an
invariant of the representation $K' \rightarrow \mathrm{GL}(V_i)$.
If $d_0$ is the minimum positive degree for such an invariant, $f$
must have degree $d \geq d_0$, because it has degree at least $d_0$ in
the variables corresponding to $V_i$ (these variables must appear 
somewhere in the polynomial because $f$ defines a smooth hypersurface).
Finally, by \Cref{semi-invariant}, we can conclude that the same $d_0$ above 
occurs as the degree of a semi-invariant for
$\bar{K}$ because $\bar{K}/\bar{K}'$ is abelian. 
(To be clear, by semi-invariant for $\bar{K} \subset \mathrm{PGL}(V_i)$,
we mean a semi-invariant of any scalar extension of $\bar{K}$ in
$\mathrm{GL}(V_i)$.)
In summary, $d \geq d_0$, where $d_0$ is the minimum positive
degree of a semi-invariant of $\bar{K}$.
\end{setup}

We will refer to the results from \Cref{elim_setup}
throughout the section and use the notation
from it, whenever there is no possibility for confusion.

The following lemma will also be useful for finding monomials excluding 
certain variables.

\begin{lemma}
\label{monomial_existence}
Let $X = \{f = 0\} \subset \mathbb{P}^{n+1}$ be a smooth hypersurface
of degree $d > 1$.  Then there is a monomial in $f$ with nonzero 
coefficient which does not involve any given set of $k$ variables, 
for any $k < \frac{n+2}{2}$.
\end{lemma}

Note that $\frac{n+2}{2} = \frac{N}{2}$, where $N$ is the dimension of the
associated vector space $V$.

\begin{proof}
Suppose every monomial in $f$ includes at least one of a given set of
$k$ variables, say $x_0,\ldots,x_{k-1}$. Then $X$ contains the 
$(n-k+1)$-plane
$\{x_0 = \cdots = x_{k-1} = 0\} \subset \mathbb{P}^{n+1}$.  But then
$n - k + 1 \leq \frac{n}{2}$ by \cite[Corollary 6.26]{3264} because
$X$ is smooth.  Rearranging
gives $k \geq \frac{n+2}{2}$.  Thus, when $k < \frac{n+2}{2}$,
there is a monomial in $f$ not containing any of the $k$ variables.
\end{proof}

\subsection{Preliminary Reductions}

This subsection is a series of lemmas that
eliminate all partitions for which exceptional
hypersurface examples do not occur. First, we prove a result about 
the degree $d$ when $\pi$ contains at least one block of size $2$.

\begin{lemma}
\label{elim2}
Suppose $G$ has a partition $\pi$ with $\mu_2 \geq 1$.  
Then 
\begin{enumerate}
    \item If $d < 12$:
$$\frac{B(\pi,d)}{|G|} \geq \min\left\{60,\left(\frac{5}{2}\right)^{\mu_2} \right\}.$$
\item If $d = 3$ and $4 \mu_2 > N$, we get further that 
$$\frac{B(\pi,d)}{|G|} \geq \min\left\{120,12\left(\frac{5}{2}\right)^{\mu_2},2 \cdot 5^{\mu_2} \right\}.$$
\end{enumerate}

\end{lemma}

\begin{proof}
First, we'll prove item (1). We consider two cases.  The first is when a central extension of $A_5$
occurs as at least one of the primitive constituents. Without loss of generality, set $\bar{H}_1 = A_5$.
Consider \Cref{elim_setup} for $i = 1$.  The kernel $\bar{K}$ 
of the composition is either trivial or all of $A_5$ since $A_5$ is simple.
If the kernel is trivial, then $H/Z(H)$ has index at least $60$ in 
$\bar{H}_1 \times \cdots \times \bar{H}_r$, so $\frac{B(\pi,d)}{|G|} \geq 60$,
so the condition in (1) is satisfied.

If $\bar{K} = A_5$, then we conclude by \Cref{elim_setup}
and \Cref{primitive2} that $d \geq 12$, contradicting the assumption 
$d < 12$ in (1).

The second case is that a central extension of $A_5$ does not appear among
the primitive constituents. This means
that for each $V_i$ of dimension $2$, $|\bar{H}_i| \leq |S_4| = 24$ by the 
classification in \Cref{primitive2}.  This reduces the maximum possible
order of $H/Z(H)$ by a factor of at least 
$$\left(\frac{60}{24}\right)^{\mu_2} = \left(\frac{5}{2}\right)^{\mu_2},$$
completing the proof.

Now we prove (2). Since $4 \mu_2 > N$,
it follows from \Cref{monomial_existence} that there are monomials in $f$ 
with nonzero coefficient and variables only from the spaces $V_i$ of 
dimension $2$. Apply \Cref{Collinsrepl} to replace $G$ with a new group 
$\hat{G}$. Writing $\hat{H}$ as in \Cref{Collinsrepl}, we know $Z(H)$ is 
a subgroup of $Z(\hat{H}) = Z(\hat{H}_1) \times \cdots \times Z(\hat{H}_n)$. 
(In this product, each 
factor $Z(\hat{H}_i)$ acts trivially on all other $V_j, j \neq i$.)
We claim that $Z(H)$ cannot equal $Z(\hat{H})$, just by virtue of the degree
being odd.  Indeed, all central extensions of $A_4, S_4$, or $S_5$
have even order centers so the matrix $-I_V \subset Z(\hat{H})$.  Hence
if $Z(H) = Z(\hat{H})$, then since $f$ is invariant under $Z(H)$, the
monomials involving only these variables from dimension $2$ summands
$V_i$ must have even degree, hence $f$ has even degree, a contradiction.  
Thus, $|Z(H)|$ falls short of the maximum $|Z(\hat{H})| \leq d^r$ by
a factor of at least $2$, because $Z(H)$ is a proper subgroup of $Z(\hat{H})$.

Now we again consider the possible primitive constituents.
We already know from the proof of (1) in the lemma
that $[H:Z(H)]$ falls short of the maximum by at least $60$ if $A_5$ is among 
the $\bar{H}_i$.  Therefore, suppose it's not.  If some primitive constituent 
is a central extension of $S_4$, without loss of generality suppose
$\bar{H}_1 = S_4$. Again consider \Cref{elim_setup} for $i = 1$. The kernel
$\bar{K}$ of the composition is then either $S_4$, $A_4$,
the Klein four-subgroup of $S_4$, or trivial.  
If it is $S_4$ or $A_4$, $d \geq 4$ by \Cref{elim_setup}
and \Cref{primitive2}, a contradiction.
Therefore, $\bar{K}$ has order at most $4$ and $[H:Z(H)]$ is less than
the product of the $|\bar{H}_i|$ by at least a factor of $24/4 = 6$, which
in turn is less than the maximum theoretical value of $[H:Z(H)]$ by 
at least $\left(\frac{5}{2}\right)^{\mu_2}$, since no $A_5$ constituents 
appear. Finally, if there are no primitive constituents 
which are $A_5$ or $S_4$, we have that the possible order $[H:Z(H)]$
is reduced by 
$$\left( \frac{60}{12} \right)^{\mu_2} = 5^{\mu_2}.$$
In summary, $[H:Z(H)]$ falls short of its maximum from 
\Cref{boundlinf} by at least $60$, at least
$6 \cdot \left(\frac{5}{2}\right)^{\mu_2}$, or at least $5^{\mu_2}$, 
depending on the case.  Since we showed $[Z(\hat{H}):Z(H)] \geq 2$ also,
we can multiply these bounds to obtain the statement in the lemma.
\end{proof}

We can use this lemma to eliminate most of the candidate partitions
in \Cref{exceptionalpart}.  For clarity, we'll divide into two corollaries,
each of which uses one of the parts of the lemma.

\begin{corollary}
\label{cor:elim2-1}
Suppose that $G = \mathrm{Lin}(f)$ has an exceptional partition which is one 
of the following:
No. 7, No. 12, No. 17, No. 19, No. 21, No. 26, No. 32, No. 34, No. 36, 
No. 38-39, No. 41-42, No. 44, No. 47, No. 49, No. 51-53, No. 55-56, 
No. 58-59, No. 61-67, No. 69-80.
Then $|G| < |B((1^N),d)|$.
\end{corollary}

\begin{proof}
In each listed case, the maximum degree listed in \Cref{exceptionalpart} is 
less than $12$ and \\
$B(\pi,3)/B((1^N),3)$ is less than 
$\min\left\{60,\left(\frac{5}{2}\right)^{\mu_2} \right\}$ 
so the result holds by \Cref{elim2}(1).
\end{proof}

\begin{corollary}
\label{cor:elim2-2}
Suppose that $G = \mathrm{Lin}(f)$ has an exceptional partition which is 
one of the following:
No. 11, No. 25, No. 33, No. 48, No. 60, No. 68.
Then $|G| < |B((1^N),d)|$.
\end{corollary}

\begin{proof}
Each of these partitions satisfies $4 \mu_2 > N$ so \Cref{elim2}(2) applies
when $d = 3$. Since $B(\pi,3)/B((1^N),3)$ is less than the minimum in the 
lemma, an exceptional example cannot occur for $d = 3$.
Thus we can assume $d \geq 4$ and we can replace the quotient 
$B(\pi,3)/B((1^N),3)$ in the table with $B(\pi,4)/B((1^N),4)$ to
apply the reasoning of \Cref{elim2}(1).  We verify
that this quotient is smaller than
$\min\left\{60,\left(\frac{5}{2}\right)^{\mu_2} \right\}$ in each case.
\end{proof}

Next we can eliminate some partitions with blocks of size $3$.

\begin{lemma}
\label{elim3}
Suppose $G$ has partition $\pi$ with $\mu_3 \geq 1$.  Then if $d < 6$:
$$\frac{B(\pi,d)}{|G|} \geq \min\left\{5,\left(\frac{15}{7}\right)^{\mu_3} \right\}.$$
\end{lemma}

\begin{proof}
Arrange the composition $V = V_1 \oplus \cdots \oplus V_r$ so that
$\dim(V_1) = 3$.  We consider two cases.  The first case is when the first primitive constituent $\hat{H}_1$
has quotient $\bar{H}_1$ equal to $A_6$ or the Hessian
group $(\Z/3)^{\oplus 2} \rtimes \mathrm{SL}_2(3)$ of order $216$.
Consider \Cref{elim_setup} for $i = 1$.  The kernel $\bar{K}$ of the 
composition there is a normal subgroup of $\bar{H}_1$.  
Since $A_6$ is simple and the largest
proper normal subgroup of $(\Z/3)^{\oplus 2} \rtimes \mathrm{SL}_2(3)$
is $(\Z/3)^{\oplus 2} \rtimes Q_8$, $\bar{K}$ either all of $\bar{H}_1$
or has order at most $216/3 = 72$.  In the latter situation, we get that
$|H/Z(H)|$ falls short of the maximum possible order by at least
$360/72 = 5$, so the condition of the lemma is satisfied.  

If the kernel is all of $\bar{H}_1$, we apply the conclusion of
\Cref{elim_setup} and the facts about semi-invariant degrees from
\Cref{primitive3} to conclude $d \geq 6$, contradicting the assumption.

Finally, if neither of these types of primitive constituents appear
among the $\hat{H}_i$, the largest possible constituent for a block of size
$3$ is a central extension of $\mathrm{PSL}_2(7)$.
In this case, the maximum
possible order of $H/Z(H)$ is reduced by a factor of at least
$$\left(\frac{360}{168}\right)^{\mu_3} = \left(\frac{15}{7}\right)^{\mu_3},$$
completing the proof of the lemma.
\end{proof}

\begin{corollary}
\label{cor:elim3}
Suppose that $G = \mathrm{Lin}(f)$ has an exceptional partition which is one 
of the following: No. 5, No. 16, No. 40, No. 50.
Then $|G| < |B((1^N),d)|$.
\end{corollary}

\begin{proof}
In each listed case, the maximum degree listed in \Cref{exceptionalpart} is
less than $6$ and \\ $B(\pi,3)/B((1^N),3)$ is less than 
$\min\left\{5,\left(\frac{15}{7}\right)^{\mu_3} \right\}$.
\end{proof}

Next, we can eliminate some partitions with blocks of size $4$.

\begin{lemma}
\label{elim4}
Suppose $G$ has partition $\pi$ with $\mu_4 \geq 1$.  Then if $d < 8$,
$$\frac{B(\pi,d)}{|G|} \geq \min \left\{720, \left( \frac{18}{5}\right)^{\mu_4} \right\}$$
\end{lemma}

\begin{proof}
As usual, there are two cases.  For the first, suppose there is a 
primitive constituent (which we can assume is $\hat{H}_1$) that
falls into case \eqref{prim4:Sp4(3)} or \eqref{prim4:H4} of 
\Cref{primitive4}, where in the second case we assume
$\hat{H}_1/Z(\hat{H}_1) = N/Z(N)$ for the full normalizer
$N$ of $H_4 \subset \mathrm{GL}_4(\C)$.  Now consider \Cref{elim_setup}
for $i = 1$.  If $\bar{H}_1 = \mathrm{PSp}_4(3)$, which is simple, the
kernel $\bar{K}$ is the whole group or trivial.
If $\bar{H}_1 = N/Z(N)$ instead, we note that any normal subgroup of 
$\bar{H}_1$ is index $2$ (corresponding to $A_6 \subset S_6$) or has index at
least $720$.
So either $[H:Z(H)]$ falls short of its maximum by at least $720$, or
$\bar{K}$ is among $\mathrm{PSp}_4(3), N/Z(N)$, or the index $2$ normal 
subgroup of $N/Z(N)$.  In this latter situation, by \Cref{elim_setup} and 
\Cref{primitive4}, $d \geq 8$, a contradiction.

The second case is that each constituent is not of the type above, so that 
by \Cref{primitive4}, $[\hat{H}_i:Z(\hat{H}_i)] \leq 7200$ for each $i$.
This means that the order of $[H:Z(H)]$ falls short of the maximum possible
by at least
$$\left(\frac{25920}{7200}\right)^{\mu_4} = \left(\frac{18}
{5}\right)^{\mu_4}.$$
This finishes the lemma.
\end{proof}

\begin{corollary}
\label{cor:elim4}
Suppose that $G = \mathrm{Lin}(f)$ has an exceptional partition which is one of the following:
No. 15, No. 24, No. 37, No. 45-46, No. 57.
Then $|G| < |B((1^N),d)|$.
\end{corollary}

\begin{proof}
In each listed case, the maximum degree listed in \Cref{exceptionalpart} is 
less than $8$ and \\
$B(\pi,3)/B((1^N),3)$ is less than $\min \left\{720, \left( \frac{18}{5}\right)^{\mu_4} \right\}$ so the result holds by \Cref{elim4}.
\end{proof}

\begin{lemma}
\label{elim6}
Suppose $G$ has partition $\pi$ with $\mu_{6} \geq 1$.  Then if $d < 6$,
$$\frac{B(\pi,d)}{|G|} \geq \min \left\{604800, \left( \frac{112}{3}\right)^{\mu_{6}} \right\}.$$
\end{lemma}

\begin{proof}
Arrange the summands so $\dim(V_{1}) = 6$. There are two cases, depending
on what $\hat{H}_1$ is.  In the first case, $\hat{H}_1$ is one of the
groups listed in \Cref{primitive6}.  In this event, we consider \Cref{elim_setup}
for $i = 1$. The possibilities for $\bar{K}$
are $\mathrm{PSU}_{4}(3).2_{2}$, its index $2$ normal subgroup $\mathrm{PSU}_{4}(3)$, 
the Janko group $J_2$, or the trivial group.  If $\bar{K}$ is trivial,
$H/Z(H)$ has index at least $|J_2| = 604800$ in the product of the
$\hat{H}_i$, giving the first bound in the lemma. 
For other possible $\bar{K}$, we conclude from \Cref{elim_setup} and 
\Cref{primitive6} that $d \geq 6$, a contradiction.

The second possibility is that the constituents above do not appear,
so by \Cref{primitive6}, the maximum possible $[H:Z(H)]$ is at most $174960$.
Therefore, the maximum possible order of $H/Z(H)$ is reduced 
from its theoretical upper bound in \Cref{boundlinf} by a factor of at least
$$\left(\frac{6531840}{174960}\right)^{\mu_{6}} = \left(\frac{112}{3}\right)^{\mu_{6}}.$$
\end{proof}

\begin{corollary}
\label{cor:elim6}
Suppose that $G = \mathrm{Lin}(f)$ has an exceptional partition which is one of the following:
No. 20, No. 28-29, No. 35, No. 43, No. 54.
Then $|G| < |B((1^N),d)|$.
\end{corollary}

\begin{proof}
In each listed case, the maximum degree listed in \Cref{exceptionalpart} is 
less than $6$ and \\
$B(\pi,3)/B((1^N),3)$ is less than $\min \left\{604800, \left( \frac{112}{3}\right)^{\mu_6} \right\}$ so the result holds by \Cref{elim6}.
\end{proof}

At this point, we only have a handful of partitions left to eliminate.
These will be treated individually.

\begin{lemma}
\label{lem:elim_sporadic}
Suppose that $G = \mathrm{Lin}(f)$ has an exceptional partition which is one of the following:
No. 3, No. 8-10, No. 14, No. 22-23, No. 27, No. 30-31.
Then $|G| < |B((1^N),d)|$.
\end{lemma}

\begin{proof}
For each partition, we assume $|G| \geq B((1^N),d)$ and derive 
a contradiction.

\noindent \textbf{Partitions No. 8: $\pi = (5)$ and No. 27: $\pi = (8)$.}

In these cases, $G$ is primitive. In order for $|G| \geq B((1^N),d)$
for some $d \geq 3$ we need $[G:Z(G)]d \geq N! d^N$ so 
$[G:Z(G)] \geq N! d^{N-1}$.  Even for $d = 3$, the right hand side
is $5! \cdot 3^4 = 9720$ or $8! \cdot 3^7 = 88179840$, precisely the 
bounds that appear in \Cref{primitive5} and \Cref{primitive8}, respectively.
But the maximum possible degree from \Cref{exceptionalpart} is $3$ 
in both cases and a semi-invariant polynomial of degree at least
$3$ must actually satisfy $d \geq 4$ for all the primitive groups listed
in \Cref{primitive5} and \Cref{primitive8}.

The next few cases make use of \Cref{monomial_existence}.

\noindent \textbf{Partition No. 3: $\pi = (2,1)$.}
By \Cref{monomial_existence}, some monomials in $f$ involve only
the variables from $V_1$ of dimension $2$, so $f(x_0,x_1,0)$ is nonzero and
invariant under the primitive representation 
$\hat{H}_1 \subset \mathrm{GL}(V_1)$.  Since $d \leq 10$ by \Cref{exceptionalpart},
$\bar{H}_1$ cannot be $A_5$.  If it is $S_4$, $d \geq 6$ by \Cref{primitive2}
but the maximum order of $|G|$ from \Cref{boundlinf}
is $24 \cdot d^2 < 6 \cdot d^3 = B((1^3),d)$ for these values of $d$.  If it is $A_4$,
$d \geq 4$ but $12 \cdot d^2 < 6 \cdot 4^3 = B((1^3),d)$ for these values of $d$.

\noindent \textbf{Partition No. 9: $\pi = (4,1)$.}
By \Cref{monomial_existence}, some monomials in $f$ involve only
the variables from $V_1$ of dimension $4$.  Hence $f(x_0,x_1,x_2,x_3,0)$
is nonzero and invariant under the primitive representation 
$\hat{H}_1 \subset \mathrm{GL}(V_1)$.  Since $B(\pi,3)/B((1^5),3) = 8$,
$[\hat{H}_1:Z(\hat{H}_1)]$ is
within a factor of $8$ of the largest possible value, $25920$.  By
\Cref{primitive4}, no group with this property has a semi-invariant of degree
at most $6$ (the max degree from the table), except for powers of the quadric
$x_0 x_3 - x_1 x_2$ in case \eqref{prim4:tensor}.  If $d = 4$ or $6$
and $f(x_0,x_1,x_2,x_3,0)$ is a power of that quadric, it would
mean a smooth threefold hypersurface of degree $4$ or $6$ would contain a 
multiple quadric as a hyperplane section, 
contradicting the Grothendieck-Lefschetz theorem, which says $\Pic(\mathbb{P}^4) \rightarrow \Pic(X)$ is an isomorphism.

\noindent \textbf{Partition No. 10: $\pi = (3,2)$.}
Let $\dim(V_1) = 3$ and $\dim(V_2) = 2$.  Consider \Cref{elim_setup} for $i = 2$. If $\bar{H}_2$ is $A_5$, the kernel $\bar{K}$
is either all of $A_5$ or trivial.
If it is trivial, $[H:Z(H)]$ is reduced by a factor of $60$ from its maximum
value, and we're done since the value $B/B_F$ in \Cref{exceptionalpart} is 6.66.  Otherwise, $d \geq 12$, which is greater than the max degree from the table, a contradiction.

Hence we've shown $|\bar{H}_2| \leq 24$, and 
\Cref{boundlinf} gives $|G| \leq 24 [\hat{H}_1:Z(\hat{H}_1)] d^2$.
Comparing this to the Fermat hypersurface, we need $[\hat{H}_1:Z(\hat{H}_1)] \geq 135$ and $d = 3$ or $[\hat{H}_1:Z(\hat{H}_1)] \geq 320$ and $d = 4$.
By \Cref{monomial_existence}, the polynomial $f(x_0,x_1,x_2,0,0)$ is
not identically zero and is invariant under $\hat{H}_1 \subset \mathrm{GL}(V_1)$.  But \Cref{primitive3} implies that there is no semi-invariant
polynomial of a primitive group in dimension $3$ satisfying one of the 
listed conditions on degree and order.

\noindent \textbf{Partition No. 14: $\pi = (4,2)$.}
Let $\dim(V_1) = 4$ and $\dim(V_2) = 2$.  Consider \Cref{elim_setup}
for $i = 2$. 
If $\bar{H}_2$ is $A_5$, the kernel $\bar{K}$ is either all of $A_5$ or trivial. If it is trivial, $[H:Z(H)]$ is reduced by a factor of $60$
from its maximum
value, and we're done by \Cref{exceptionalpart}.  Otherwise, $d \geq 12$, which is greater than the max degree from the table. 

Hence for an exceptional example we need $|\bar{H}_2| \leq 24$, and 
\Cref{boundlinf} gives $|G| \leq 24 [\hat{H}_1:Z(\hat{H}_1)] d^2$.  Comparing
this to the Fermat hypersurface, we need either $d = 4,5$ and 
$[\hat{H}_1:Z(\hat{H}_1)] \geq 7680$ or $d =3$
and $[\hat{H}_1:Z(\hat{H}_1)] \geq 2430$.  The polynomial
$f(x_0,x_1,x_2,x_3,0,0)$, which is nonzero by \Cref{monomial_existence}, 
is an invariant of the primitive representation $\hat{H}_1 \subset \mathrm{GL}(V_1)$ and would have to satisfy the order and degree conditions above. \Cref{primitive4} rules out the existence of such a polynomial.

\noindent \textbf{Partition No. 22: $\pi = (4,3)$.}
Let $\dim(V_1) = 4$ and $\dim(V_2) = 3$.  Consider \Cref{elim_setup}
for $i = 2$.  Suppose first that $\bar{H}_2$ is $A_6$ or the Hessian
group.  The possible kernels $\bar{K}$ are $\bar{H}_2$ itself, the index
$3$ normal subgroup of the Hessian, or a group of order at most $18$.
(This is because the next largest normal subgroup of the Hessian group
has order $18$.)
In the latter scenario, the order of $H/Z(H)$ falls short of its maximum
by a factor of at least $360/18 = 20$, which is greater than the value
$7.62$ from \Cref{exceptionalpart}.  But if $\bar{K}$ is one of the other
groups listed, $d \geq 6$, which is greater than the maximum of $4$ from the 
table.  Hence the largest possible $\bar{H}_2$ is $\mathrm{PSL}_2(7)$.

Thus \Cref{boundlinf} gives $|G| \leq 168 [\hat{H}_1:Z(\hat{H}_1)] d^2$.
Comparing this to the Fermat hypersurface, we get $d = 3$ and 
$[\hat{H}_1:Z(\hat{H}_1)] \geq 7290$. Using \Cref{monomial_existence}, the polynomial $f(x_{0}, x_{1}, x_{2}, x_{3}, 0, 0, 0)$ is nonzero and is an invariant of the primitive representation $\hat{H}_1 \subset \mathrm{GL}(V_1)$.  By \Cref{primitive4}, no group satisfies the necessary condition
on order and degree.

\noindent \textbf{Partition No. 23: $\pi = (4,2,1)$.}
Let $\dim(V_1) = 4, \dim(V_2) = 2$, and $\dim(V_3) = 1$. 
Consider the composition from \Cref{elim_setup} with $i = 2$.
If $\bar{H}_2$ is $A_5$, the kernel $\bar{K}$ is all of $A_5$ or trivial. If it is trivial, $[H:Z(H)]$ is reduced by a factor of $60$ from its maximum
value, and we're done by \Cref{exceptionalpart}.  Otherwise, $d \geq 12$, which is greater than the max degree from the table.

Hence for an exceptional example we need $|\bar{H}_2| \leq 24$, and 
\Cref{boundlinf} gives $|G| \leq 24 [\hat{H}_1:Z(\hat{H}_1)] d^3$.
Comparing this to the Fermat hypersurface, this forces $d = 3$ and
$[\hat{H}_1:Z(\hat{H}_1)] \geq 17010$.  The polynomial
$f(x_0,x_1,x_2,x_3,0,0,0)$, which is nonzero by \Cref{monomial_existence}, 
is an invariant of the primitive representation $\hat{H}_1 \subset \mathrm{GL}(V_1)$; this polynomial and $\hat{H}_1$ would have to satisfy the degree and order conditions above. \Cref{primitive4} rules out the existence of such a polynomial.

\noindent \textbf{Partition No. 30: $\pi = (4^2)$.}
In this case, we have two primitive constituents $\hat{H}_1$ and 
$\hat{H}_2$, both primitive subgroups of $\mathrm{GL}_4(\C)$.
Suppose one constituent, say $\hat{H}_1$ without loss of generality,
satisfies $\bar{H}_1 = \mathrm{PSp}_4(3)$ or $\bar{H}_1 = N/Z(N)$, for
the group $N$ from \Cref{primitive4} \eqref{prim4:H4}.
Then the kernel $\bar{K}$ of the composition in \Cref{elim_setup} for $i = 1$
is a normal subgroup of $\bar{H}_1$, so that it is
all of $\mathrm{PSp}_4(3)$ (in case \eqref{prim4:Sp4(3)}), either $N/Z(N)$ or its index $2$ subgroup (in case \eqref{prim4:H4}), or has index at least $720$ in $\bar{H}_1$.  When the index is at least $720$, $[H:Z(H)]$
falls short of the upper bound by at least that number, so we're done
by \Cref{exceptionalpart}.  

In the remaining cases, we use \Cref{elim_setup} to conclude $d \geq 8$, which is larger than the maximum degree value from the table.  Hence $[\hat{H}_i:Z(\hat{H}_i)] \leq 7200$ for each $i$. 
On the other hand, \Cref{boundlinf} gives
\begin{equation}
\label{No. 30}
   |G| \leq 2 [\hat{H}_1:Z(\hat{H}_1)][\hat{H}_2:Z(\hat{H}_2)] d^2. 
\end{equation}
Comparison to $B((1^8),d)$ yields $d = 3$ and 
$$[\hat{H}_1:Z(\hat{H}_1)][\hat{H}_2:Z(\hat{H}_2)] \geq 14696640.$$
The factor of $2$ in front of \eqref{No. 30} comes from the possibility
of a symmetry between the two factors, but this can only occur if
$\hat{H}_1 \cong \hat{H}_2$.
Consulting the list in \Cref{primitive4} gives that the only remaining possibility is that $\hat{H}_i/Z(\hat{H}_i) = Q/Z(Q)$ for each $i$,
where $Q$ has an index $2$ subgroup $2.A_5 \times 2.A_5/Z$ and $Z$ is a central subgroup of order $2$.
(Here $Q$ is the largest group from \Cref{primitive4}\eqref{prim4:tensor}).

The last thing to do is see that $d = 3$ cannot occur for this choice
of primitive constituents.
Returning to \Cref{elim_setup} again for $i = 1$, $\bar{K}$ is a 
normal subgroup of $Q/Z(Q)$ which is either the whole subgroup, a subgroup of 
index $2$, or trivial. 
If it is trivial, $[H:Z(H)]$ falls short of the upper
bound by at least $7200$ and
we're done.  In the other two cases, the conclusion of \Cref{elim_setup}
gives that every term of $f$ has even degree in the variables $x_0, x_1,x_2,x_3$ corresponding to $V_1$. But we can repeat this whole argument
with \Cref{elim_setup} again for $i = 2$ and get that
every term of $f$ has even degree in $x_4,x_5,x_6,x_7$.  Hence $f$ has even degree, and $d \neq 3$, a contradiction.

\noindent \textbf{Partition No. 31: $\pi = (4,2^2)$.}
Let $\dim(V_1) = 4$ and $\dim(V_2) = \dim(V_3) = 2$.  We'll work
by finding stronger bounds on the orders of each primitive constituent.
Suppose that at least one of the dimension $2$ constituents, say $\hat{H}_2$,
is an extension of $A_5$.  Then the group $\bar{K}$ in \Cref{elim_setup} for
$i = 2$ is trivial or all of $A_5$.  If it is trivial, the upper bound on 
$[G:Z(G)]$ decreases by $60$ and we're done by \Cref{exceptionalpart}.
If it is all of $A_5$, $d \geq 12$, contradicting the max
degree from the table.  Hence each of the two-dimensional
constituents satisfies $[\hat{H}_i:Z(\hat{H}_i)] \leq 24$.  Now \Cref{boundlinf} gives
$$|G| \leq 2 \cdot 24^2 [\hat{H}_1:Z(\hat{H}_1)]d^3.$$
Comparing to the Fermat hypersurface gives that $d = 3$ and $[\hat{H}_1:Z(\hat{H}_1)] \geq 8505$.
Hence $\hat{H}_1$ must be as in \eqref{prim4:Sp4(3)} or \eqref{prim4:H4} in
\Cref{primitive4} (in the latter case, $\hat{H}_1:Z(\hat{H}_1) \cong N/Z(N)$ for the entire normalizer $N$, in particular).
Now we consider \Cref{elim_setup} for $i = 1$, so that $\bar{K}$ is a normal subgroup of $\bar{H}_1$.  This must either be
all of $\bar{H}_1$, have index $2$, or else its index is so large that
$|H/Z(H)|$ is reduced by a factor of at least $720$, and we'd be done.
In the former cases, we'd have $d \geq 8$, contradicting $d = 3$.
\end{proof}

\subsection{The Exceptional Examples}

Putting together Corollaries \ref{cor:elim2-1}, \ref{cor:elim2-2},
\ref{cor:elim3}, \ref{cor:elim4}, \ref{cor:elim6}, and 
\Cref{lem:elim_sporadic}, the only
remaining exceptional partitions to consider are No. 1, No. 2, No. 4, No. 6,
No. 13, and No. 18.  We won't consider No. 1, the partition $\pi = (2)$,
because this would correspond to a zero-dimensional hypersurface $X$.

The remaining cases are partitions of $N = 3,4$, or $6$,
so in particular we've finished proving that no $X$ with
corresponding partition $\pi \neq (1^N)$ can match or exceed
the Fermat unless $\dim(X) = n = 1,2,4$.
Each of the remaining partitions will actually yield hypersurfaces
with larger automorphism group than the Fermat hypersurface in certain
degrees.  In this section, we'll find which degrees can occur and 
prove that the examples from \Cref{sect:mainresults} are unique, up
to isomorphism.

\begin{proof}[Proof of \Cref{mainthmfull}:]
We'll consider each of the remaining partitions in turn.
We begin each case by assuming $|G| \geq B((1^N),d)$.

\noindent \textbf{Partition No. 2: $\pi = (3)$.}

In this case, the group $G = \mathrm{Lin}(f)$ is primitive,
so it fits into one of the cases of \Cref{primitive3}. Further,
$|G| = [G:Z(G)] d$.  In contrast, the linear automorphism
group in $\mathrm{GL}_3(\C)$ of the Fermat equation of degree $d$
has order $|G_F| = 6d^3$.

Hence any exceptional example would have 
$[G:Z(G)] d \geq 6d^3$, or $[G:Z(G)] \geq 6d^2$.  In particular,
$[G:Z(G)] \geq 6 \cdot 3^2 = 54$.
There are four possibilities, following \Cref{primitive3}.

If $G/Z(G) = H/Z(H)$ for $H = 3.A_6$, then $[G:Z(G)] = 360$
and hence $d \leq 7$.  The unique semi-invariant polynomial in this range is
the Wiman sextic, and conversely, the Wiman sextic does define a smooth
plane curve.  Hence $X$ is isomorphic to Example
\ref{ex:(1,6)}.  Since $\mathrm{Aut}(X) \cong A_6$, $|\mathrm{Aut}(X)| = 360$.

If $G/Z(G) = H/Z(H)$ for a subgroup $H$ of the Hessian group with $[H:Z(H)] 
\geq 54$, $d \leq 6$ since $[H:Z(H)] \leq 216$. But the smallest degree
semi-invariants also have degree $6$ in this case, so $d = 6$.  The largest 
possible $\mathrm{Aut}(X)$ is the full Hessian group of order $216$, which
exactly equals the order of the automorphism group of the Fermat sextic.
Hence $X$ is isomorphic to Example \ref{ex:(1,6)-2}.

If $G/Z(G) = H/Z(H)$ for $H = \mathrm{PSL}_2(7)$, then $[G:Z(G)] = 168$,
so $d \leq 5$.  By \Cref{primitive3}, the unique semi-invariant polynomial
in this range is the Klein quartic, and conversely, the Klein quartic
does define a smooth plane curve with $|\mathrm{Aut}(X)| = 168$.  Hence
$X$ is isomorphic to Example \ref{ex:(1,4)}.

If $G/Z(G) = H/Z(H)$ for $H = 2.A_5$, then $[G:Z(G)] = 60$ so $d = 3$.
But there is no semi-invariant of degree $3$.

\noindent \textbf{Partition No. 4: $\pi = (4)$.}

In this case, the group $G = \mathrm{Lin}(f)$ is again primitive,
so $|G| = [G:Z(G)] d$.  In contrast, the linear automorphism group
in $\mathrm{GL}_4(\C)$ of the Fermat equation of degree $d$ has 
order $|G_F| = 24 d^4$.

Hence any exceptional example would have 
$[G:Z(G)] d \geq 24d^4$, or $[G:Z(G)] \geq 24d^3$.  Since $d \geq 3$,
in particular it's true that $[G:Z(G)] \geq 648$, exactly the 
condition in \Cref{primitive4}.  We need only check these cases.

We see that in cases \eqref{prim4:Sp4(3)}, \eqref{prim4:2A7}, and
\eqref{prim4:2S6} of \Cref{primitive4}, that $[G:Z(G)] < 24d^3$, 
where $d$ is the minimum degree of a semi-invariant.  In case \eqref{prim4:tensor},
when $1152 < [G:Z(G)] \leq 7200$, $d \leq 6$ and there are possible
semi-invariants in degrees $4$ or $6$, but they are all powers of the quadric
$x_0 x_3 - x_1 x_2$ and do not define a \textit{smooth} surface.
If $[G:Z(G)] \leq 1152$, $d < 4$ so there are no available semi-invariants
to define a surface.

Lastly, in case \eqref{prim4:H4}, $[G:Z(G)] \leq 11520$ so $d \leq 7$.
We therefore can't get any semi-invariants of the required degree unless 
$G/Z(G) \cong H/Z(H)$ for $H \subset N$ the subgroup corresponding to the 
exceptional embedding of $S_5 \subset S_6$ (or its index $2$ subgroup).  
Either of these subgroups satisfies $[H:Z(H)] \leq 11520/6 = 1920$, so 
$d \leq 4$, and $d \leq 3$ for the index $2$ subgroup.  By the conditions
on semi-invariants from \Cref{primitive4}, we must have $d = 4$, and
the equation $f$ is the unique quartic from that proposition.  Therefore,
$X$ is isomorphic to Example \ref{ex:(2,4)}.

\noindent \textbf{Partition No. 6: $\pi = (2^2)$.}

In this case, the group $G = \mathrm{Lin}(f)$ is not primitive, but
instead has two primitive constituents $\hat{H}_1, \hat{H}_2$
of dimension $2$. The only possibilities for these constituents are central
extensions of $A_4$, $S_4$, or $A_5$. By \Cref{boundlinf} we have
\begin{equation}
\label{No.6}
  |G| \leq 2 [\hat{H}_1:Z(\hat{H}_1)][\hat{H}_2:Z(\hat{H}_2)] d^2.  
\end{equation}
On the other hand, we need $|G| \geq 24 d^4$, so 
$[\hat{H}_1:Z(\hat{H}_1)][\hat{H}_2:Z(\hat{H}_2)] \geq 12 d^2$.
Suppose $\bar{H}_1 = A_5$.  Using \Cref{elim_setup} for $i = 1$, the kernel
$\bar{K}$ is trivial or all of $A_5$.  If it is trivial,
$|G| \leq 2 [\hat{H}_2:Z(\hat{H}_2)] d^2$ and we'd get $d < 3$.  Thus, the kernel is all of $A_5$, $d \geq 12$,
and $f$ is an invariant polynomial of a central extension of $A_5$ in the 
first two variables $x_0,x_1$. 
Only the choice $\bar{H}_1 = \bar{H}_2 = A_5$ gives a large
enough group for this value of $d$, and each term of $f$ must symmetrically
be an invariant polynomial in the variables $x_2,x_3$ of the second summand 
as well.  Finally, we know $d \leq 17$ from \Cref{exceptionalpart}, so the 
only possibility is that $f$ is the sum of the unique degree $12$ invariant
in both pairs of variables.  
Therefore, $X$ is isomorphic to Example \ref{ex:(2,12)}.

If neither constituent is an extension of $A_5$, consider the
case  $\bar{H}_1 = \bar{H}_2 = S_4$.  Using \eqref{No.6} again
gives $d \leq 6$.  We examine the possibilities for the subgroup
$H/Z(H) \subset \bar{H}_1 \times \bar{H}_2$.  If it is the whole group,
by \Cref{elim_setup} applied to either projection, $d = 6$ and 
$f$ must be an invariant polynomial in each pair of variables, forcing it
to be the sum of the unique degree $6$ invariant in both pairs of variables.
Therefore, $X$ is isomorphic to Example \ref{ex:(2,6)}.

If $H/Z(H)$ is a proper subgroup of the product, it is either not
symmetric under the interchange of factors (so we remove the $2$ in the front
of the bound in \eqref{No.6}) or has index at least $4$. The same reduction
in order happens when we change the constituents: if $\bar{H}_1 = A_4$
and $\bar{H}_2 = S_4$, we must also remove the factor of $2$, and if 
$\bar{H}_1 = \bar{H}_2 = A_4$, then $[\hat{H}_1:Z(\hat{H}_1)]
[\hat{H}_2:Z(\hat{H}_2)] = 144$.  In all of the cases from this paragraph,
$$|G| \leq 288 Z(H) \leq 288 d^2$$
so $d = 3$. But the oddness of the degree means $Z(H)$ is a proper subgroup
of $Z(\hat{H}) = Z(\hat{H}_1) \times Z(\hat{H}_2)$ (using \Cref{Collinsrepl}), so
in fact $|G| \leq 144 d^2$.  This is too small to compete with the Fermat
for any degree $d \geq 3$.

\noindent \textbf{Partition No. 13: $\pi = (6)$.}

In this case, the group $G = \mathrm{Lin}(f)$ is primitive and
an exceptional example must satisfy 
$|G| = [G:Z(G)] d \geq 720 d^6$, where the latter is the order of the
linear automorphism group of the Fermat equation in dimension $6$.  This
gives $[G:Z(G)] \geq 720 d^5$.  Plugging in $d = 3$ gives 
$[G:Z(G)] \geq 174960$, exactly the condition from \Cref{primitive6}. 
So $G/Z(G) \cong H/Z(H)$ for one of the two groups from this proposition.
In case \eqref{prim6:J2}, the smallest degree possible is $12$, which
is much too large, leaving only \eqref{prim6:U4(3)}.  The
requirement $[G:Z(G)] \geq 720 d^5$ for $[G:Z(G)] = 6531840$ gives $d \leq 6$,
so the only possible equation is the unique one stated in \Cref{primitive6}.
Therefore, $X$ is isomorphic to Example \ref{ex:(4,6)}.

\noindent \textbf{Partition No. 18: $\pi = (2^3)$.}

In this case, the group $G = \mathrm{Lin}(f)$ is not primitive, but
instead has three primitive constituents $\hat{H}_1, \hat{H}_2, \hat{H}_3$
of dimension $2$.  By \Cref{boundlinf} we have
\begin{equation}
\label{No.18}
    |G| \leq 6 [\hat{H}_1:Z(\hat{H}_1)][\hat{H}_2:Z(\hat{H}_2)][\hat{H}_3:Z(\hat{H}_3)] d^3.
\end{equation}
In contrast, the Fermat equation has automorphism group of 
order $|G_F| =  720 d^6$.  Even setting the max value of $60$ for
each $[\hat{H}_1:Z(\hat{H}_1)]$, we must have $d \leq 12$.
First, suppose at least one primitive
constituent $\hat{H}_1$ is a central extension of $A_5$.  Then applying
\Cref{elim_setup} for $i = 1$, the kernel $\bar{K}$ is all of $A_5 = \bar{H}_1$ or trivial.  If it is trivial, this reduces
the maximum possible order of $H/Z(H)$ by $60$ and the modified \eqref{No.18}
gives $d \leq 3$.  But if $d = 3$, the degree being odd forces
$|Z(H)| \leq \frac{d^3}{2}$, 
and $6 \cdot 60^2 \cdot d^3/2$ is smaller than $720 d^6$, even for $d = 3$.

So $\bar{K}$ is all of $A_5$, and $d = 12$ is forced by \Cref{elim_setup}.
The only choice of primitive constituents consistent with this bound is 
$\bar{H}_1 = \bar{H}_2 = \bar{H}_3 = A_5$.  Since we can apply the same
argument to each of the three factors, the only possible equation is the sum
of the unique degree $12$ invariant polynomial in each pair of variables.  
Therefore, $X$ is isomorphic to Example \ref{ex:(4,12)}.

If no primitive constituent is an extension of $A_5$,
suppose $\bar{H}_1 = S_4$.  
Unless the other two constituents are also extensions of $S_4$, 
the group $G/H$ cannot be a transitive subgroup of $S_3$ and so the 
coefficient $6$ in \eqref{No.18} is reduced to at most $2$.  In this case,
$d < 3$, a contradiction. 
(We also get $d < 3$ if all three constituents are extensions of $A_4$.) 
So $\bar{H}_1 = \bar{H}_2 = \bar{H}_3 = S_4$ and the inequality
in \eqref{No.18} gives $d \leq 4$.

Considering possible subgroups $H/Z(H) \subset \bar{H}_1 \times \bar{H}_2 \times \bar{H}_3$, either $H/Z(H)$ is this whole group, forcing $d \geq 6$ (a contradiction) or $H/Z(H)$ is a proper subgroup.  This proper subgroup is
either not fully symmetric under the permutation of factors, or it has index
at least $8$.  Either way, the bound on $|G|$ drops to the point that $d < 3$, completing the proof.

\noindent \textbf{Uniqueness for other $(n,d)$.}

We've now dealt with all exceptional cases.
Suppose once again that the hypersurface $X_d = \{f = 0\} \subset \mathbb{P}^{n+1}$ has 
automorphism group of order at least
that of the Fermat hypersurface. In all other dimensions $n = N-2$
and degrees $d$ besides the seven pairs
listed in \Cref{mainthmfull}, we've proven 
that the group $G = \mathrm{Lin}(f)$ has order $|G| < B((1^N),d)$
unless the partition $\pi$ of $G$ equals $(1^N)$.
So if $(n,d)$ is not one of the
seven special pairs, then the associated $G = \mathrm{Lin}(f)$ has partition
$\pi = (1^N)$.  Further, $|G|$ must achieve the upper bound $B((1^{N}),d)$.
In light of the analysis of this section, we still have $\pi = (1^N)$
and $|G| = B((1^{N}),d)$ even if $(n,d)$ is one of the seven special pairs,
unless $X$ is isomorphic to the eight examples
listed in \Cref{mainthmfull}.  The last claim to check is that
when $\pi = (1^N)$ and $|G| = B((1^{N}),d)$, $X$ must in fact
be isomorphic to the Fermat hypersurface.

Because $|G|$ already equals $B((1^N),d)$, applying Collins' 
replacement theorem (\Cref{Collinsrepl}) to produce a new group 
$\hat{G}$ outputs a group of the same order.  
Following the notation of that theorem, since $H \subset \hat{H}$, we 
must have $H = \hat{H} = \hat{H}_1 \times \cdots \times \hat{H}_N$, 
and $G/H = S_N$ is the full symmetric group.  Each $\hat{H}_i$ is a 
finite subgroup of $\mathrm{GL}(V_i) \cong \C^*$ and acts trivially 
on all other $V_j$. This means that each term of $f$ must have degree
a multiple of $|\hat{H}_i|$ in each variable $x_i$.  The product of 
the $|\hat{H}_i|$ is $d^N$ and $f$ has degree $d$, so the only possibility 
is $\hat{H}_i = \mu_d$ for each $i$.  Thus, $f$ is a linear combination 
of the monomials $x_i^d$.  Finally, $f$ is invariant under the symmetric
group permuting the subspaces $V_i$, so $f$ must in fact be the Fermat polynomial.
\end{proof}

\end{document}